\documentclass[11pt,a4paper]{article}

\renewcommand{\baselinestretch}{1}

\usepackage{fullpage}
\usepackage{titlesec}
\titleformat*{\section}{\large\bfseries}

\usepackage{latexsym,amssymb,amsthm,amsmath,amsfonts,graphicx,booktabs} 
\usepackage[round]{natbib}
\usepackage[linesnumbered,lined,ruled,longend]{algorithm2e}

\newcommand{\diag}{\text{diag}}
\newcommand{\Var}{\text{var}\,}
\newcommand{\Cov}{\text{cov}\,}
\newcommand{\X}{\mathcal{X}}
\newcommand{\Tc}{\mathcal{T}}
\newcommand{\T}{\mathrm{\scriptscriptstyle T}}
\newcommand{\tr}{\text{tr}}
\newcommand{\supp}{\text{supp}}

\newtheorem{theorem}{Theorem}

\newtheorem{remark}{Remark}

\usepackage{collab}
\collabAuthor{wm}{red!85}{Werner}
\collabAuthor{mh}{magenta!85}{Markus}
\collabAuthor{ap}{blue!85}{Andrej}

\usepackage{xcolor}
\definecolor{rev1}{rgb}{0,0,0}

\begin{document}

\begin{center}
\renewcommand{\baselinestretch}{1}
{ \large
{\bf A convex approach to optimum design of experiments with correlated observations}
}

\vspace{0.5cm}

Andrej~P\'azman$^a$,
Markus~Hainy$^b$ and
Werner~G.~M\"uller$^b$\footnote{Corresponding author: Department of Applied Statistics, Johannes Kepler University Linz, Altenberger Stra{\ss}e 69, 4040 Linz, Austria; \texttt{werner.mueller@jku.at}} \\[3ex]
$^a$ Comenius University Bratislava, Slovakia \\[1ex]
$^b$ Johannes Kepler University Linz, Austria 
\end{center}

\renewcommand{\baselinestretch}{1}
%\small 

\noindent Abstract: Optimal design of experiments for correlated processes is an increasingly relevant and active research topic. Present methods have restricted possibilities to judge their quality. To fill this gap, we complement the virtual noise approach by a convex formulation leading to an equivalence theorem comparable to the uncorrelated case and to an algorithm giving an upper performance bound against which alternative design methods can be judged. Moreover, a method for generating exact designs follows naturally. We exclusively consider estimation problems on a finite design space with a fixed number of elements. A comparison on some classical examples from the literature as well as a real application is provided.     
  \\

%\normalsize
%\vspace{-0.5cm}

%\newpage
%\vspace{-0.5cm}
\section{Introduction and motivation}\label{sec:intro}
%\vspace{-1cm}

Due to its importance in spatio-temporal monitoring (cf. \citealp{mateu_spatio-temporal_2012-1})
and particularly computer simulation (cf. \citealp{rasmussen_gaussian_2005}), 
experimental design for regression models with correlated errors has gained increasing interest. {The theory is well developed for classical (non-)linear regression with uncorrelated errors (cf. \citealp{morris_design_2010}), where the corresponding optimum design approach has been initiated from Kiefer's concept of design measures (cf. \citealp{kiefer_optimum_1959}) and has emerged into textbook status (cf. 
	\citealp{atkinson_optimum_2007-1}, \citealp{goos_optimal_2011}, etc.). However, for the correlated setting the literature is only scattered.}

The meanwhile classic approach by \cite{sacks_designs_1966} relied heavily on asymptotic considerations, making it valid and useful for only limited situations. In contrast, the algorithm devised by \cite{brimkulov_numerical_1980} is purely heuristic but applicable to and surprisingly efficient for a great variety of problems. 
The proposal of \cite{fedorov_design_1996} to transform the problem into a random coefficient model allows for embedding it into standard convex design theory but requires elaborate tuning to achieve a given design size. More recently, in a remarkable series of papers started with \cite{zhigljavsky_new_2010}, Dette, Pepelyshev and Zhigljavsky built upon and extended much of the discussions and material exposed in the pioneering monograph by \cite{nather_effective_1985}, essentially concentrating on the ordinary least squares estimator. Latest, more promising approaches derived from the best linear unbiased estimator of the continuous-time model were proposed in \cite{dette_optimal_2016-1} and \cite{dette_new_2017}.

An entirely different suggestion redefining the role of design measures was given in \cite{pazman_new_1998} and fully developed in \cite{muller_measures_2003}. Therein a virtual noise was introduced to influence the behavior of the designs and the corresponding measure was reflecting the amount of signal suppression. The resulting procedures were quite effective, but their broader application was hampered by the downside that the method was nonconvex and thus no Kiefer-Wolfowitz-type equivalence theorem was available. In the present paper we fill this gap by asserting convexity through a new and different variant of virtual noise and consequently provide a corresponding equivalence theorem.

We will show the effectiveness of this novel modification on classical examples from the literature and compare to the above mentioned alternative techniques. 

\section{Problem description and theory}
\subsection{The setup} \label{subsec:setup}

In the presentation we will constrict ourselves on linear regression models with the obvious generalizations to 
the asymptotic approach to the nonlinear case as outlined in \cite{pronzato_design_2014}, so we assume our observations to be generated from the model
\begin{equation} \label{setup}
y\left( x\right) =f^\T\left( x\right) \theta + \varepsilon\left( x\right)
;\;\;\;x\in \mathcal{X},
\end{equation}
where $\mathcal{X}=\{x_1,\dots,x_N\}$ is a finite discrete design space  containing $N$ points and $\theta $ is the unknown vector of parameters with dim$\left(
\theta \right) =p$.  What is further conventionally assumed is that the random error of the observation
performed at $x$ has zero expectation, i.e. $E\left\{ \varepsilon\left( x\right) \right\} =0$, and is correlated, defined through 
the $N\times N$ matrix $C$ with elements
\[
C_{ij}=\Cov\left\{\varepsilon\left( x_i\right) , \varepsilon\left(
x_j\right) \right\} ;\quad \quad \,\,i,j=1,\ldots,N, 
\]
which is supposed to be known. 
The Gaussian process literature additionally assumes normality, which leads to the well-known kriging equations for prediction (cf. \citealp{rasmussen_gaussian_2005}). 

Let us now first consider an exact unreplicated $n$-point design $\Tc \subset \X$, where $n$ is the required number of
observations.
The information matrix of $\theta $ in the model (\ref{setup}) and for the design $\Tc$ is inversely related to the (asymptotic) variance-covariance matrix of the best linear unbiased estimator $\hat{\theta}$. As is well known, this information matrix is  
$
M_\Tc =F(\Tc)^\T C ^{-1}(\Tc)F(\Tc)$,
where $C(\Tc)$ is the $n\times n$ submatrix of $C$ corresponding to points of the set $\Tc$ and $F(\Tc)$ is an $n\times p$ matrix, 
$F_{i,j}(\Tc)=f_j\left( x_i\right)$, with $x_i \in \Tc$ and $j=1,\ldots,p$.

The situation is quite different when the observations are uncorrelated, so that for any $x \in \X$ observations $y (x)$ can be replicated. Traditionally, following \cite{kiefer_optimum_1959},
we then define a design as a probability measure {$\xi$} on $\X$, $\xi (x)$ being interpreted
as to be proportional to the number of replicated observations at the design point
$x$.

Typically criteria functions used for 
uncorrelated
observations, i.e. the case of a diagonal $C$, can be expressed as concave increasing functions of the
information matrix. For example, for D-optimality we have 
$
\Phi \left( M_\xi\right) =\left\{ \det \left( M_\xi\right) \right\} ^{1/p} 
\mbox{ or } 
\Phi \left( M_\xi\right) =\log \det \left( M_\xi\right), 
$
with \linebreak
$
M_{\xi} = \sum_{x\in \mathcal{X}} f(x) f^\T(x) \, \xi(x). \label{FIMuncorr}
$
{
	It is the usual purpose of optimum design to find the design 
	\begin{equation*} \label{xistar}
	\xi^*=\arg\max_{\xi} \Phi \left( M_\xi\right).
	\end{equation*}
	The maximum of $\Phi \left( M_\xi
	\right) $ over the set of all designs $\Xi $ can be obtained by convex methods. In case
	of strict concavity of $\Phi \left( M_\xi \right) $, there is only one maximum value of $M$.
}

It is fully justified to use in the correlated setup of model (\ref{setup}) the
same optimality criteria functions $\Phi(M)$ as for the uncorrelated case, but now
$M_\Tc$ should be used instead of $M_\xi$ (\citealp{pazman_criteria_2007}). In the method presented in this paper, we
take great advantage of the fact that all statistically meaningful criteria functions $\Phi(M)$ are concave functions on the set of all symmetric positive semidefinite
matrices $M$, cf.~\cite{Pukelsheim_2006}. 
Thus, we will be able to employ efficient convex optimization methods for the evaluation of approximate designs. 

Note that this is usually not the case in the correlated setup, i.e. for a general $C$, where typically sophisticated non-convex optimization techniques, such as simulated annealing or particle swarms, have to be applied.

\subsection{Some techniques from the literature}

Such non-convex methods are for instance required for the mathematically elaborate design approaches developed by Dette, Zhigljavsky and collaborators, most recently in \cite{dette_optimal_2016-1} and \cite{dette_new_2017}. There they generate their results by moving from the discretely indexed model (\ref{setup}) to its continuously indexed counterpart and back, thus ending up with the eventual need for discretisation and approximation of the best linear unbiased estimator.
In \cite{dette_optimal_2016-1}, they {make use of the fact that the optimally-weighted signed least squares estimator has the same variance as the best linear unbiased estimator. They derive the weight function of the signed least squares estimator for the continuous model and regard this as the design measure.} In
\cite{dette_new_2017}, they propose an approximate criterion based on a discrete approximation of the best linear unbiased estimator for the continuous model.  
The approximate criterion needs to be optimized numerically in the same way as the original criterion, so there is no major conceptual difference when it comes to optimization with the exception of the former being computationally more efficient if $p<<n$. 
Therefore, in the following we shall not consider these approximate criteria, as \cite{dette_new_2017} have already investigated how good/bad this approximation works.
%%%%%%%%%%%%%%%%%%%%%%%%%%%%%%%%%%%%%%%%%%%%%%%%%%%%%%%%%%%%

Limitations of the methods of \cite{dette_optimal_2016-1} and \cite{dette_new_2017} are that they are constrained to a limited form of covariance kernels and that it may be quite difficult to generalize to higher-dimensional design spaces.

In contrast, the suggestion by \cite{fedorov_design_1996} is to approximate directly the error component through Mercer's covariance expansion of degree $q$, i.e. $\sum_{l=1}^q\gamma_{l}\varphi_l(x)$, where the $\varphi_l(x)$ and
$\lambda_l$ are the eigenfunctions and eigenvalues, respectively, of
the covariance kernel. 
This leads to model (\ref{setup}) being approximated by a mixed (random and fixed) coefficient regression model and allows for embedding into standard convex design theory, albeit involving a potentially high number of $q$ additional regressors and corresponding inflation of the information matrix plus some arbitrariness in determining an appropriate $q$.   Incidentally it was shown in \cite{pazman_note_2010} that for a finite setup the method can be embedded into {techniques similar to the one} described in the next subsection with very specific variances of the virtual noise components. 

The most straightforward and pragmatic approach, however,  is the greedy algorithm for exact design generation first proposed by \cite{brimkulov_numerical_1980}, in the modification and interpretation given in \cite{fedorov_design_1996}.
There the so-called sensitivity function of the classical one-point correction algorithms for D-optimality in the uncorrelated setup is used with expressions simply replaced by their analogues from the correlated case. That is, at each step $r$ we augment the point
$$
x_r = \arg \max_{\X}  f^\T(x)M^{-1}_{\Tc_{r-1}}f(x) / \sigma^2_\varepsilon(x)
$$ 
to the current design {$\Tc_{r-1} \subset \Tc$}, with $f$, $M$, and the error variance $\sigma^2_\varepsilon$ replaced accordingly, see the appendix for details. Furthermore, the appendix contains an adaption for A-optimality, which will be utilized in some of the examples. Note that despite its compellingly simple form, this algorithm has little theoretical basis other than that the independent case can be seen as a particular limit instance of the general setup.  

A good overview of all these approaches with an application in spatio-temporal sensor placement, albeit mainly from the perspective  of using the ordinary least squares estimator, can most recently be found in \cite{ucinski_d-optimal_2020}.

\subsection{Virtual noise}

In this paper, we will use a technique that has first been proposed in \cite{pazman_new_1998}, which is very different from the above. Despite a certain similarity to the standard kriging setup, in that we consider two error components, one correlated and the other uncorrelated, our additional uncorrelated error  component 
	however is not considered as an observational error but as a regulatory device without physical meaning, termed virtual noise. The purpose of this noise is to suppress the signal whenever the design measure is small relative to the maximum value of the design measure.
To make this operational, the variance of this virtual noise needs to be assumed being of a specific form. Unfortunately the specifications given in  \cite{pazman_new_1998} as well as in \cite{muller_measures_2003} both yielded a nonconvex solution. 
In this paper, we intend to fill this gap by providing convexity and consequently an equivalence theorem that allows for quick checks of whether particular designs optimize the design criterion as well as an algorithm to find the optimum.

First, similarly as in \cite{muller_measures_2003}, 
define a convex set of restricted probability measures, also termed design measures or simply designs, 
\[
\Xi =\left\{ \xi :\sum_{x\in \mathcal{X}}\xi \left( x\right) =1, \: \forall
_{x\in \mathcal{X}} \: 0~\leq~\xi \left( x\right) \leq 1/n\right\}, 
\]
on the design space $\mathcal{X}$, which is supposed to be finite. The integer number $n$ represents the desired number of points
of support of an optimum exact design $\xi^*$. 
Now instead of (\ref{setup}),
consider the potentially perturbed variables 
\begin{equation} \label{perturbed}
y\left( x\right) =f^\T\left( x\right) \theta +\varepsilon \left( x\right)
+w_{\xi} \left( x\right) ;\;\;\;x\in \mathcal{X},
\end{equation}
where the variance of the supplementary 'virtual
noise' $w_{\xi}$, independent of $\varepsilon$, is fixed as
\[
\Var\left\{ w_{\xi}\left( x\right) \right\} = \sigma^2_{\xi}(x) =\kappa \frac{1/n-\xi \left( x\right) }{
	\xi \left( x\right)}, 
\]
whereas
$
\Cov \left\{ w_{\xi}\left( x\right), w_{\xi}\left( x'\right) \right\} = 0 \text{ when } x \neq x'. 
$
The positive number $\kappa$ is a tuning parameter which must be chosen within the bounds given by Theorem~\ref{concav} and  should be taken as large as possible to emphasize the influence of the virtual noise. The experience obtained in the examples below is that a smaller choice of $\kappa$ diminishes mainly the speed of the linear programming algorithm from Section~\ref{subsec:design_algorithm}.

That means, similarly as in \cite{muller_measures_2003}, if $\xi \left( x\right) =0$, i.e. $\sigma^2_\xi = \infty$, there is no observation at the point 
$x$, and if $\xi \left( x\right) =1/n$, i.e. $\sigma^2_\xi = 0$, the observation at the point $x$ is
not disturbed at all by the virtual noise.

An important notion is the concept of an exact $k$-point design, for which in our context we require a slightly different definition as usual. Here, it is a design $\xi \in \Xi$ for which the signal is not disturbed by the virtual noise in the $k$ points $x_1,\ldots,x_k \in \mathcal{X}$, i.e.~$\xi(x_i) = 1/n$, and the signal is totally suppressed in all other points, i.e.~$\xi(x) = 0 \; \: \forall \: x \in \mathcal{X} \backslash \{x_1,\ldots,x_k\}$. It follows trivially from the definition of $\Xi$ that all $n$-point exact designs are in $\Xi$, but $\Xi$ contains no $k$-point exact designs with $k \neq n$. The set $\Xi$ has also another property: it is the smallest closed convex set of probability measures on $\mathcal{X}$ which contains all $n$-point exact designs. The primary aim of the statistician is to consider exact $n$-point designs, so the set $\Xi$ is the smallest possible convex extension of the set of such exact designs.

In the proofs we shall also require the set
	\begin{equation*}
	\Xi_+ = \left\{ \xi : \sum_{x \in \mathcal{X}} \xi(x) = 1, \: \forall_{x \in \mathcal{X}} \: 0 < \xi(x) \leq 1/n  \right\} = \left\{ \xi \in \Xi : \: \supp(\xi) = \mathcal{X} \right\}.
	\end{equation*}
	Its advantage is that $\Var\left\{ w_{\xi}\left( x\right) \right\} < \infty$ for every $x \in \mathcal{X}$, $\xi \in \Xi_+$. One can easily see that $\Xi$ is the closure of $\Xi_+$.

Denote by $W\left( \xi \right) $ the $N\times N$ diagonal matrix 
\begin{eqnarray*}
	W_{i,j}\left( \xi \right) &=&\kappa \frac{1/n-\xi \left( x_i\right) }{\xi
		\left( x_i\right) };  \mbox{ if }\,\,\,i=j \\
	&=&0; \qquad \qquad \quad \: \: \mbox{ if }\,\,i\neq j.
\end{eqnarray*}

When $\xi \in \Xi_+$, i.e. $\supp(\xi) = \mathcal{X}$,
	the information matrix of $\theta$ in model (\ref{perturbed}) is given by 
\begin{equation} \label{FIMvirtual}
M( \xi ) =F^\T\left\{ C+W\left( \xi \right) \right\} ^{-1}F,
\end{equation}
where $C$ is the $N \times N$ covariance matrix defined after (\ref{setup}), and $F$ is the $N \times p$
matrix $F_{i,j} = f_j (x_i),\: x_i \in \X,\: i = 1,\ldots,N,\: j = 1,\ldots, p$. Evidently, $M(\xi)$ is nonsingular if $F$ has full rank $p$.
	
	On the other hand, when $\xi \in \Xi\backslash\Xi_+$, we define 
	\begin{equation} \label{FIMvirtual2}
		M( \xi ) =F'^\T\left\{ C'+W'\left( \xi \right) \right\} ^{-1}F',
	\end{equation}
	where $C'$, $W'(\xi)$ are submatrices of $C$, $W(\xi)$ restricted to the support $\supp(\xi)$, and similarly are the rows of $F'$. 
	An important property is the continuity of $M(\xi$) on the whole set $\Xi$, i.e. $\xi_n(x) \to \xi(x)$ for every $x \in \mathcal{X}$ implies $M(\xi_n) \to M(\xi)$, see Lemma A2 in the appendix.

	Another interesting property of $M(\xi)$ is that in case of uncorrelated observations with constant variances, $C=\sigma^2I$, with $\kappa=\sigma^2$ we obtain $M(\xi)=M_\xi$, the information matrix from Kiefer's design theory. However, for general $C$  the design measure can no more be interpreted as reflecting the replications of observations. In any case we can still use it to obtain heuristically exact $n$-point designs approaching the quality of an optimum exact design as will be demonstrated in the examples section. 

We now have the following 
\begin{theorem}
	\label{concav}
	If  $
	\kappa < \lambda _{\min }\left( C\right), 
	$
	the minimal eigenvalue of the matrix $C,$ and if $\Phi
	\left( M\right) $ is any optimality criterion expressed as a
	concave, increasing, and continuous function of the matrix $M$, then the mapping  
	\[
	\xi \in \Xi \rightarrow \Phi \left\{ M( \xi ) \right\}
	\]
	is concave as well, with $M( \xi )$ defined
	in $\left(\ref{FIMvirtual}\right)$ and $\left(\ref{FIMvirtual2}\right)$.
\end{theorem}

The proof of  the theorem is given in the appendix. Its consequences are that the maximum of $\Phi(M)$ over the set $\Xi $ can be obtained by convex methods, and in case
of a strict concavity of $\Phi(M)$, there is only one optimal $M$.

Unfortunately however, the gradient method so popular in uncorrelated observations cannot
be used here because one-point design measures do not enter into the set $
\Xi $. Neither can we use the usual equivalence theorem, which requires
the use of one-point designs as well. However, the methods of 
modified linear programming can be employed, see the next subsection, and also a certain form of the
equivalence theorem can be provided, which is given in Section~\ref{subsec:equi_theorem}.

\subsection{A design algorithm} \label{subsec:design_algorithm}

The method used in this subsection originated as the cutting-plane method of \cite{kelley_cutting-plane_1960}, was then adapted to experimental design in  Chapter 9 of \cite{pronzato_design_2014} and applied algorithmically to classical experimental design in \cite{burclova_optimal_2016}.

For computational reasons we shall restrict our attention in this section to the set
	\begin{equation*}
	\Xi^{\epsilon} = \left\{ \xi \in \Xi : \: \forall_{x \in \mathcal{X}} \: \xi(x) \geq \epsilon\right\}
	\end{equation*}
	for some small $\epsilon > 0$. Evidently,
	\[ \Xi_+ = \bigcup_{\epsilon > 0} \Xi^{\epsilon}. \]
	Since the closure of $\Xi_+$ is the set $\Xi$ and since according to Lemma~A2 in the appendix the criterion function $\Phi \left\{ M( \xi ) \right\}$ is continuous on the whole of $\Xi$, we can approach the value of $\max_{\xi \in \Xi} \Phi \left\{ M( \xi ) \right\}$ arbitrarily close by computing $\Phi \left\{ M( \xi^*_{\epsilon} ) \right\}$, where $\xi^*_{\epsilon} \in \arg \max_{\xi \in \Xi^{\epsilon}} \Phi \left\{ M( \xi ) \right\}$.

If $\Phi ( M ) $ is a global criterion, like $D-$ or $A-$optimality,  which has a gradient 
$
\nabla _M\Phi ( M ), 
$
then for any $\xi$ in $\Xi_+ =\{\xi \in \Xi: \supp(\xi)=\X\}$  we have the derivative 
\[
\frac{\partial \Phi \left\{ M( \xi ) \right\} }{\partial
	\xi ( \bar{x}) }=\tr\left[ \nabla _M\Phi \left\{ M (
\xi ) \right\} \,\,\frac{\partial M( \xi ) }{
	\,\,\,\partial \xi( \bar{x}) }\right],
\]
where $\nabla _M\Phi \left\{ M (\xi ) \right\}$ denotes $\nabla _M\Phi ( M )$ at $M=M(\xi)$.
The linear Taylor formula for $\xi \rightarrow \Phi \left\{ M
( \xi ) \right\} $ at the point $\mu $ gives
\begin{equation}
\Phi \left\{ M( \mu ) \right\} +\sum_{\bar{x}\in \mathcal{X}
}\tr\left[ \nabla _M\Phi \left\{ M ( \mu ) \right\} \,\,
\frac{\partial M( \mu )}{\,\,\,\partial \mu ( \bar{x
	}) }\right] \left\{ \xi ( \bar{x} ) -\mu ( \bar{x}
) \right\}, \label{linear_Taylor_formula}
\end{equation}
which is linear in $\xi $, so due to the concavity, we have for any $\xi \in \Xi^{\epsilon}$
\begin{equation}
\Phi \left\{ M\left( \xi \right) \right\} =\min_{\mu \in \Xi^{\epsilon} }\left(
\Phi \left\{ M( \mu ) \right\} +\sum_{\bar{x}\in \mathcal{X}
}\tr\left[ \nabla _M\Phi \left\{ M ( \mu ) \right\} \,\,
\frac{\partial M( \mu ) }{\,\,\,\partial \mu ( \bar{x
	}) }\right] \left\{ \xi ( \bar{x} ) -\mu ( \bar{x}
) \right\} \right) \label{linear_Taylor_formula_min}
\end{equation}
and $\arg \max_{\xi} \Phi\{M(\xi)\}$ is the optimal design on $\Xi^{\epsilon}$. It follows that for a differentiable positive criterion function $\Phi$
the optimal design problem is an ``infinite-dimensional linear programming problem'', which
can be solved iteratively.

To see the details, let us take a look at the particular specification for D-optimality.
Take $\Phi( M ) = \left\{\det ( M ) \right\}^{1/p} $. 
Then for any nonsingular $M$
\[
\nabla _M\Phi (M) = M^{-1} \cdot p^{-1} \left\{ \det ( M ) \right\}^{1/p}.
\]

According to (\ref{FIMvirtual}), supposing that $F$ has full rank we have that $M^{-1}( \mu ) $
exists for any $\mu \in \Xi_+$, and we also have
\begin{equation*}\label{deltaM}
\frac{\partial M( \mu ) }{\,\,\,\partial \mu ( \bar{x
	}) }=\frac \kappa n F^\T H^{-1}( \mu ) \Gamma \left\{ \mu
^{-2}( \bar{x}) \right\} H^{-1}( \mu ) F, 
\end{equation*}
with
\begin{equation*}\label{H}
H( \mu ) =\left[ \left( C-\kappa I\right) +\frac \kappa
n \, \diag\left\{ \mu ^{-1}( \cdot ) \right\} \right], 
\end{equation*}
where $\diag\left\{ \mu ^{-1}( \cdot ) \right\}$ is an $N\times N$
diagonal matrix with $\mu ^{-1}( x ) $ as diagonal elements, and 
$
\Gamma \left\{ \mu^{-2}( \bar{x}) \right\}
$ is an $N\times N$ matrix having zero elements everywhere up to the $\bar{x}$th
diagonal position, where we have $\mu^{-2}( \bar{x}) $.
Thus, the expression \eqref{linear_Taylor_formula} becomes
\begin{eqnarray*}
	&& 
	\left[ \det \left\{ F^\T H^{-1}( \mu ) F\right\} \right]^{1/p} \times \\ && \times 
	\left[ 1 
	+\frac{\kappa}{n p} \sum_{\bar{x}\in \mathcal{X}} \left\{H^{-1}( \mu )\right\}_{\bar x,.}
	F  \left\{ F^\T H^{-1}( \mu ) F\right\}^{-1} F^\T 
	\left\{H^{-1}( \mu )\right\}_{.,\bar x}
	\frac{\xi ( \bar{x} ) -\mu ( \bar{x} ) }{\mu
		^2 ( \bar{x} ) } \right].
\end{eqnarray*}

Returning to the general case, we can now reformulate the design problem into a linear program. 
Any concave and positive criterion 
\[
\xi \in \Xi^{\epsilon} \rightarrow \Phi \left\{ M( \xi ) \right\}
\]
can be written in a form 
\[
\Phi \left\{ M ( \xi ) \right\} =\min_{\mu \in \Xi^{\epsilon} }\left\{ a (
\mu ) +\sum_{i=1}^N b_i ( \mu ) \, \xi ( x_i )
\right\} .
\]
Here, $a( \mu ) $ and $b_i ( \mu ) $ can be obtained
either from the linear Taylor formula \eqref{linear_Taylor_formula_min}, or, in case that there is no gradient, from subgradients, or using some algebraic tricks as in
\cite{burclova_optimal_2016}.

At the first step of the linear program, we set a finite set $\Xi _1\subset \Xi^{\epsilon} $, which will
be increased in the later steps.
Hence, at the $k$th step we start with a set $\Xi _k$. Consider the following
finite set of constraints, which are linear in the variables $t\in \mathbb{R}$ and $
\xi ( x );\: x\in \mathcal{X}$, 
\begin{eqnarray*}
	t &\geq&0, \qquad
	t \leq a ( \mu ) +\sum_{i=1}^N b_i ( \mu ) \, \xi (
	x_i ) \text{,\thinspace \thinspace \thinspace \thinspace \thinspace }
	\mu \in \Xi _k; \\
	\xi ( x_i ) & \geq & \epsilon;
	\qquad \xi ( x_i ) \leq \frac 1n; \quad \quad i=1,\dots,N, \qquad 
	\sum_{i=1}^N \xi ( x_i ) =1,
\end{eqnarray*}
and consider a standard linear program which maximizes $t$ under these
constraints. Denote $t^{\left( k\right) }$ and $\xi ^{\left( k\right)
} $ those values of $t$ and of $\xi $ where the maximum of $t$ is attained.
Define $\Xi _{k+1}=\Xi _k\cup \left\{ \xi^{(k)} \right\} $ and repeat the linear program
with $\Xi _{k+1}$ instead of $\Xi _k$.

To make this operational, we also require a stopping rule. We have that
\begin{eqnarray*}
	t^{\left( k\right) } &=&\max_{\xi \in \Xi^{\epsilon} }\min_{\mu \in \Xi _k}\left\{
	a ( \mu ) +\sum_{i=1}^N b_i ( \mu ) \, \xi ( x_i )
	\right\} 
	\\
	&\geq &\max_{\xi \in \Xi^{\epsilon} }\min_{\mu \in \Xi^{\epsilon}}\left\{
	a ( \mu ) +\sum_{i=1}^N b_i ( \mu ) \, \xi ( x_i )
	\right\} \\
	&=& \max_{\xi \in \Xi^{\epsilon} } \, \Phi \left\{ M ( \xi ) \right\} 
	\geq \max_{\xi \in \Xi _k}\Phi \left\{ M ( \xi ) \right\}.
\end{eqnarray*}

Hence, if we have a small number $\delta >0$, such that
$
t^{\left( k\right) }-\max_{\xi \in \Xi _k}\Phi \left\{ M ( \xi ) \right\} < \delta,
$
then 
\[
\max_{\xi \in \Xi _k}\Phi \left\{ M ( \xi ) \right\} > \max_{\xi
	\in \Xi^{\epsilon} } \, \Phi \left\{ M ( \xi ) \right\} -\delta 
\]
and 
$
\arg \max_{\xi \in \Xi _k}\Phi \left\{ M ( \xi ) \right\} 
$
is the nearly optimum design.
Evidently, $\max_{\xi \in \Xi _k}\Phi \left\{ M ( \xi ) \right\} $
can be simply computed iteratively,
\[
\max_{\xi \in \Xi _k}\Phi \left\{ M ( \xi ) \right\} =\max \left[
\max_{\xi \in \Xi _{k-1}}\Phi \left\{ M ( \xi ) \right\} ,\,\,\Phi
\left\{ M\left( \xi ^{\left( k\right) }\right) \right\} \right] .
\]

\subsection{Equivalence theorem}\label{subsec:equi_theorem}

In this section, we shall consider only global optimality criteria having a gradient $\nabla_M \Phi(M)$ such as the D- and A-criterion. We shall use an abbreviated notation writing
\begin{equation*}
	h\left( x,{\xi}\right) = \left\{ T\left( {\xi}\right) \right\}
_{x,.}G\left({\xi}\right) \left\{ T^\T\left({\xi}\right) \right\}
_{.,x},
\end{equation*}
where
\begin{eqnarray*}
	G ( \xi ) &=& F \left[ \nabla _M\Phi \left\{ M (\xi) \right\}
	\right] F^\T, \\
	T ( \xi ) &=&\left[ \left( C-\kappa I\right) \,\,\diag\left\{ \xi
	( \cdot ) \right\} +\frac \kappa nI\right] ^{-1}.
\end{eqnarray*}

\begin{theorem}\label{equivt1}
	Let $\Phi$ be a global, concave and nonnegative optimality criterion, zero on singular matrices, and having a gradient.
	\textit{
		A design measure } $\bar{\xi} \in \Xi$ with a nonsingular $M\left(\bar{\xi}\right)$  is $\Phi -$\textit{optimal
		(within our theory) if and only if for every n-tuple }$z_1,\ldots,z_n$ \textit{
		of points from }$\mathcal{X}$ \textit{we have}
	
	\begin{eqnarray} \label{equiv1}
	&&\frac 1n\sum_{i=1}^n h(z_i,\bar\xi)
	\leq \sum_{x\in \mathcal{X}} \bar\xi(x) h(x,\bar\xi).
	\end{eqnarray}

	Moreover, if $\bar{\xi}$ satisfies the inequality
	\begin{equation}\label{equiv1_delta}
	\frac 1n\sum_{i=1}^n h(z_i,\bar\xi)
	\leq \sum_{x\in \mathcal{X}} \bar\xi(x) h(x,\bar\xi)+\delta
	\end{equation}
	for some $\delta > 0$ and any $z_1,\dots, z_n \in \X$, we have that
	\begin{equation}\label{equiv1_phi}
	\Phi\{M(\bar\xi)\} \ge \max_{\xi \in \Xi} \Phi\{M(\xi)\}- \delta.
	\end{equation}
\end{theorem}

The full proof is given in the appendix. 

\begin{remark}

The numerical exploitation of this theorem is straightforward:

	For a fixed design $\bar{\xi} \in \Xi$ and for any $x\in \mathcal{X}$ compute the
	number 
	$
	h \! \left( x,\bar{\xi}\right)$
	and compute the sum 
	$
	d \! \left( \bar{\xi}\right) =n\sum_{x\in \mathcal{X}}\bar{\xi} ( x ) \,
	h \! \left( x,\bar{\xi}\right)$. Order the points of the set $\mathcal{X}$ so that 
	$
	h \! \left( x_1,\bar{\xi}\right) \geq
	\cdots \geq h \! \left( x_N,\bar{\xi}\right)$.
Then \eqref{equiv1} is equivalent to the inequality
$$
\sum_{i=1}^n h \! \left( x_i,\bar{\xi}\right) \leq d \! \left( \bar{\xi}\right) ,
$$
and \eqref{equiv1_delta} is equivalent to the inequality
$$
\sum_{i=1}^n h \! \left( x_i,\bar{\xi}\right) \leq d \! \left( \bar{\xi}\right) + n \delta.
$$
\end{remark}

Note that the possibility to formulate an equivalence theorem due to concavity means that we can provide a rather close upper bound for the performance of all exact $n$-point designs. 
Thus we can calibrate the performance of all existing methods by this bound, as will be done in the next section.

\section{Examples} \label{sec:examples}
\vskip-10pt
\subsection{General considerations} 

To check our ideas, we have computed designs for numerous examples from the literature, particularly those already presented and compared in
\cite{glatzer_comparison_1999}. In most of these examples, all compared methods did reasonably well,
so for illustration we have just selected four, which we consider representative and for which we give details below.

The ultimate goal in all our examples is to find an exact design and to judge its quality. As for all measure-based design approaches, we therefore require a method to convert the calculated measure to a discrete set of points. \cite{dette_optimal_2016-1} typically use the endpoints plus quantiles of their measure in their one-dimensional examples. We will thus for comparison purposes use the same method for our approach, but there is no simple way to extend this to higher dimensions. Instead, in the spirit of the random design strategies defined in \cite{waite_minimax_2020}, one can perform random sampling according to the measure multiple times and select the sampled design with the best criterion value. We will denote this approach by R-VN (random sampling virtual noise). In our experience, this strategy has led to even more efficient designs in all cases, and we apply this approach using 100 samples whenever the design dimension is greater than one. We furthermore provide the results for the one-dimensional cases in the appendix, which also contains some additional examples.

In the following examples, we will compare our method based on the optimal measure for the virtual noise representation to the method from  \cite{dette_optimal_2016-1} as well as to the variant of the algorithm of \cite{brimkulov_numerical_1980} put forward by \cite{fedorov_design_1996} based on approximate sensitivity functions, henceforth denoted by BKSF.
An algorithmic description of the latter can be found in the appendix. We refrain from reporting results from the original algorithm which is tantamount to an exchange method based on direct comparisons of criterion values, as most of the time they gave very similar results. 
For comparison purposes, we also provide the median performance out of 100 randomly chosen exact designs of the respective size and denote this approach by R-UNIF.

As \cite{dette_optimal_2016-1}, we only consider examples with one design variable. For all methods except the one from \cite{dette_optimal_2016-1}, the design space to search over was discretized to a $101$-point design grid. For small $n$ we performed an exhaustive search of all $n$-point combinations over the design grid to obtain the true optimum exact design. We report those points in the tables as well, where we abbreviate this approach by EXS (exhaustive search). We consider two different ways how to obtain the exact $n$-point designs for our method. The first is to select the quantiles $1/(n+1), \ldots, n/(n+1)$ of the design grid with respect to the design measure.
The second way is to always include the endpoints of the design grid and select the quantiles $1/(n-1),\ldots,(n-2)/(n-1)$ of the remaining mass over the design grid. We denote these two methods by Q-VN and Q-VN+EP (quantile virtual noise plus endpoints), respectively. 
The second approach is similar to the approach of \cite{dette_optimal_2016-1}, which we will denote by Q-DPZ+EP. The design measures \cite{dette_optimal_2016-1} consider usually have discrete masses at the endpoints but are continuous elsewhere, so \cite{dette_optimal_2016-1} always include the endpoints and select the quantiles for the design points in between with respect to the continuous measure.

From Theorem~\ref{equivt1} it follows that $\Phi\left(M_{\Tc}\right)$ for any $n$-point exact design $\Tc$ is bounded from above by $\Phi\left\{M \! \left(\bar{\xi}\right)\right\}$ for the $\Phi$-optimal design measure $\bar{\xi}$, where $M(\xi)$ is given by (\ref{FIMvirtual}). Therefore, we report the efficiencies with respect to $\Phi\left\{M \! \left(\bar{\xi}\right)\right\}$, where $\bar{\xi}$ is found by applying the linear programming algorithm.
Let $\Xi_k$ be the finite set of design measures determining the set of linear constraints in step $k$ of the linear programming algorithm and let 
	$$ t^{\left( k\right) } = \max_{\xi \in \Xi^{\epsilon} }\min_{\mu \in \Xi_k}\left\{
	a ( \mu ) +\sum_{i=1}^N b_i ( \mu ) \, \xi ( x_i )
	\right\},$$ 
	see Section~\ref{subsec:design_algorithm}. Denote $\Phi^{(k)} = \max_{\xi \in \Xi_k} \Phi\left\{M(\xi)\right\}$. In our programs, we stop the linear programming algorithm if $\frac{t^{(k)} - \Phi^{(k)}}{\Phi^{(k)}} \leq 10^{-4}.$ 
	We set $\epsilon = 10^{-6}$ in all our examples to avoid numerical problems. To solve the linear programming problems, we use the package \texttt{lpSolve} \citep{lpSolve} from the software \texttt{R} \citep{RSoftware}. In all our examples we also choose $\kappa$ to be rounded down from $\lambda_{min}(C)$ to two significant digits.

In all the figures we present, the left panel displays the discrete measure obtained by running the linear programming algorithm for the virtual noise representation for a particular $n$. 
	The right panel shows the D-efficiencies of exact designs obtained by various methods up to $n = 20$. The methods depicted are: Q-VN (solid red line with squares), Q-VN+EP (dashed green line with large dots), Q-DPZ+EP (long-dashed blue line with triangles), BKSF  (dotted black line with diamonds), and reporting the median efficiency for R-UNIF (long-short-dashed grey line with small dots).

\subsection{Example 1: a classic one-parameter model}\label{sec:example1}

This example has originally been considered by \cite{sacks_designs_1966}. \cite{dette_optimal_2016-1} use it to illustrate the efficiency of their method. It is a one-parameter model given by
\begin{eqnarray*}
	f(x) & = & 1 + 0.5 \sin(2 \pi x), \qquad     x  \in  [1,2], \\
	\Cov\{\varepsilon(x),\varepsilon( x')\} & = & \begin{cases}
		x^2 x' \qquad x \leq x' \\
		x (x')^2 \quad x > x', 
	\end{cases} \quad
	\lambda_{\min}(C)  =  0.00276. 
\end{eqnarray*}

\begin{figure}[hbtp!] 
	\vskip-20pt
	\centering
	\begin{tabular}{cc}
		\includegraphics[width=0.45\textwidth]{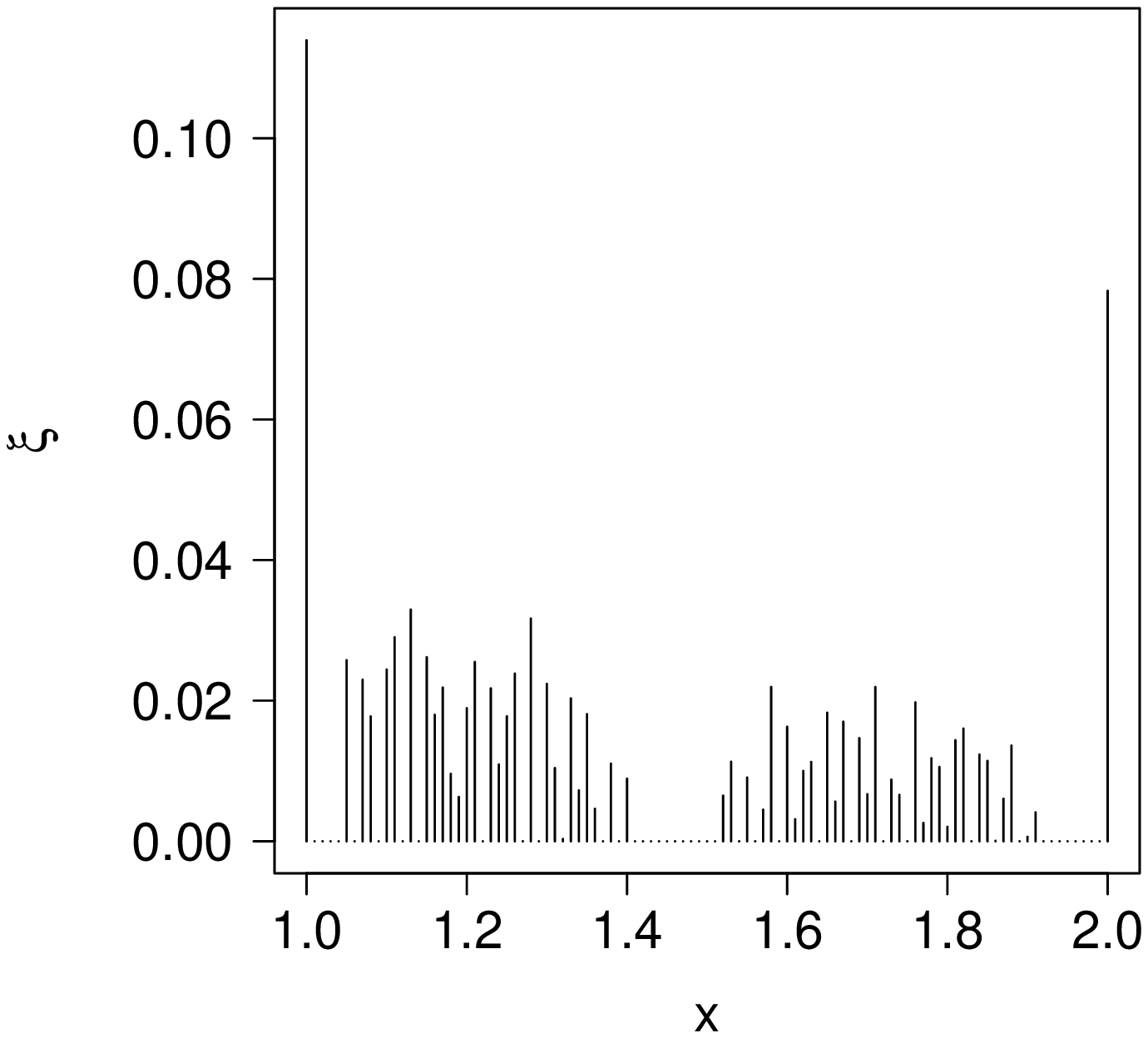} & 
		\includegraphics[width=0.45\textwidth]{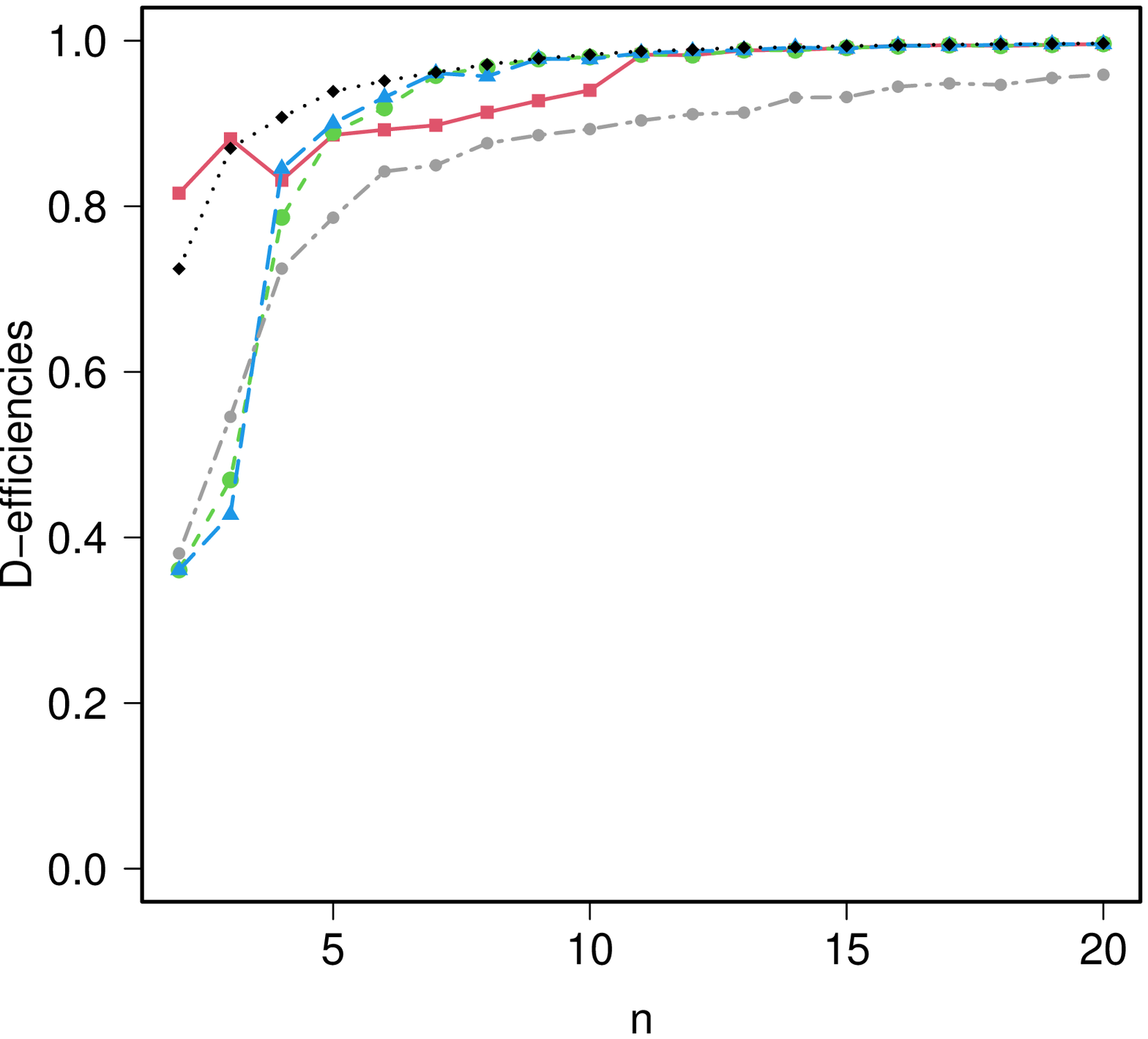}
	\end{tabular}
	\vskip-10pt
	\caption{Our measure (left panel) and efficiencies (based on our bound) versus sample size (right panel) for Example~1.
		\label{fig:measures_example1}}
\end{figure}

For the D-criterion and $n = 4$, the selected design points and D-efficiencies for the various methods are given in Table~\ref{tab:example1}. 
While the algorithm using the approximate sensitivity function from \cite{fedorov_design_1996} comes close to the result of the exhaustive search and to $90\%$ of the upper bound, the measure-based methods are considerably worse. However, when we increase the sample size (see Fig.~\ref{fig:measures_example1}), their performance improves and slowly approaches the bound.

% latex table generated in R 4.0.2 by xtable 1.8-4 package
% Mon Oct 05 16:21:37 2020
\begin{table}[hbtp!]
	\begin{center}
	\caption{Optimal designs and D-efficiencies for Example~1}
	\begin{tabular}{lccccc}
			& $x_1$ & $x_2$ & $x_3$ & $x_4$ & D-eff \\ 
			Q-VN & 1.10 & 1.23 & 1.40 & 1.76 & 0.8316 \\ 
			Q-VN+EP & 1.00 & 1.21 & 1.58 & 2.00 & 0.7865 \\ 
			Q-DPZ+EP & 1.00 & 1.28 & 1.69 & 2.00 & 0.8455 \\ 
			R-UNIF  &  &  &  &  & 0.6955 \\ 
			BKSF & 1.19 & 1.67 & 1.79 & 2.00 & 0.9075 \\  
			EXS & 1.22 & 1.66 & 1.79 & 2.00 & 0.9158 \\  
	\end{tabular}
	\end{center}
	\label{tab:example1}
\end{table}

Note that in contrast to the method by \cite{dette_optimal_2016-1} our method does not require the existence of the continuous best linear unbiased estimator. To illustrate this, we modify the example by using the kernel $\Cov\{\varepsilon(x),\varepsilon( x')\}  = \min(x,x')^2 \left\{ 3 \max(x,x') - \min(x,x') \right\} / 6$ instead, which is the once continuously differentiable kernel for the integrated Brownian motion error process. According to \cite{dette_blue_2019}, the continuous best linear unbiased estimator therefore has to incorporate information from the first derivatives of the error process to be estimable. Since this kernel is now much smoother, the minimum eigenvalue for our design grid is rather small, namely $\lambda_{\min}(C) = 2.0854 \cdot 10^{-8}$. However, we did not observe any numerical issues due to that fact.

\begin{table}[hbtp!]
	\begin{center}
	\caption{Optimal designs and D-efficiencies for the modified Example~1}
	\begin{tabular}{lccccc}
			& $x_1$ & $x_2$ & $x_3$ & $x_4$ & D-eff \\ 
			Q-VN & 1.00 & 1.01 & 1.39 & 1.53 & 0.4933 \\ 
			Q-VN+EP & 1.00 & 1.22 & 1.53 & 2.00 & 0.7329 \\ 
			R-UNIF &  &  &  &  & 0.4887 \\ 
			BKSF & 1.00 & 1.39 & 1.80 & 2.00 & 0.8042 \\ 
			EXS & 1.00 & 1.23 & 1.75 & 2.00 & 0.9715 \\  
	\end{tabular}
	\end{center}
	\label{tab:example4}
\end{table}

\begin{figure}[hbtp!]
	\centering
	\begin{tabular}{cc}
		\includegraphics[width=0.45\textwidth]{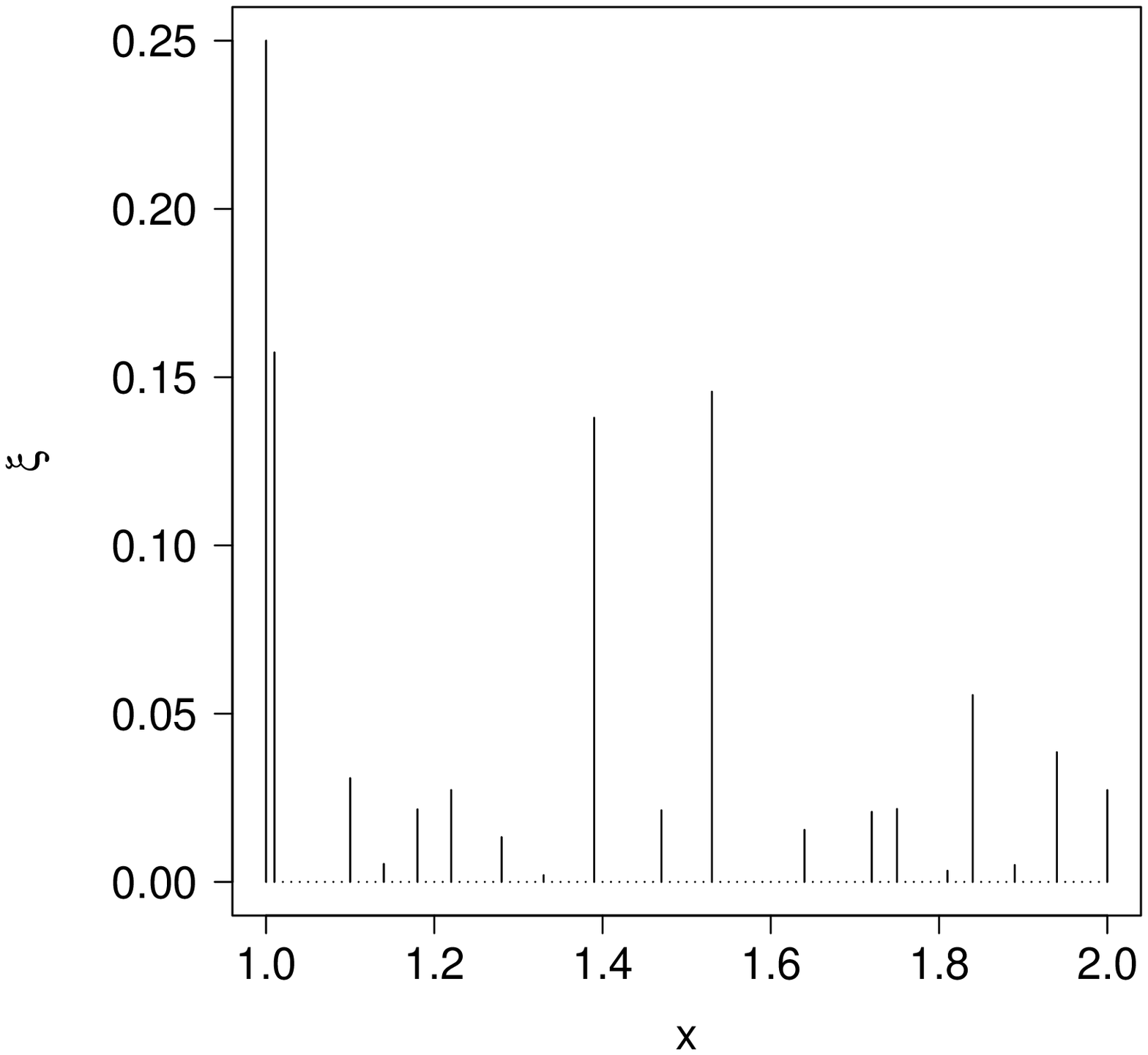} & 
		\includegraphics[width=0.45\textwidth]{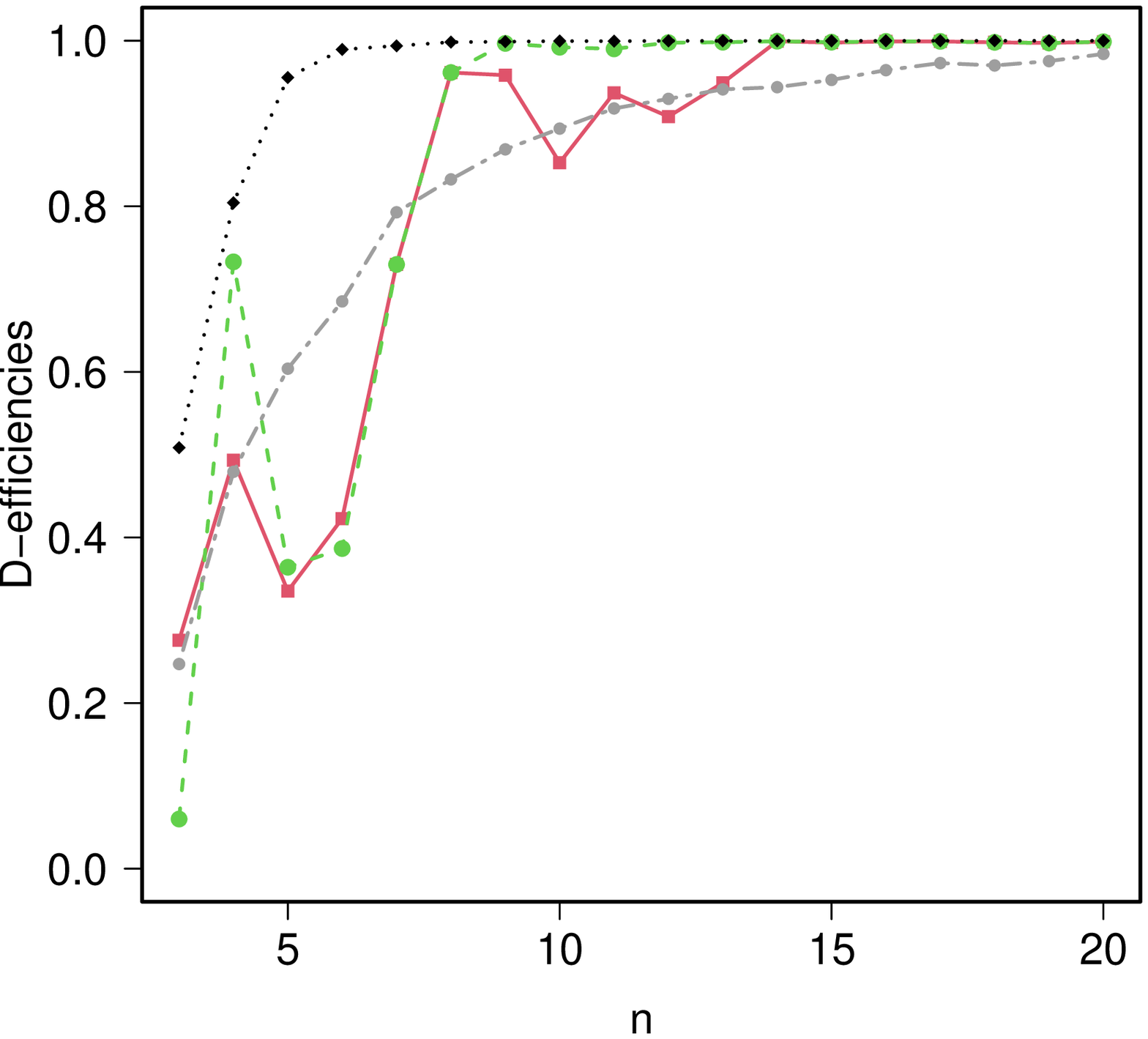}
	\end{tabular}
	\caption{Our measure (left panel) and efficiencies versus sample size (right panel) for the modified Example~1. \label{fig:measures_example4}}
\end{figure}

	Our algorithm puts a high amount of mass at the two lowest points in the design grid, see Figure~\ref{fig:measures_example4}. This can be interpreted as the algorithm trying to obtain information about the first derivative at the lower bound.

%%%%%%%%%%%%%%%%%%%%%%%%%%%%%%%%%%%%%%%%%%%%%%%%%%%%%%%%%%%%%%%%%%%%%%%%%

\subsection{Example 2: a multiparameter case}\label{sec:example2}

The next example is taken from Section~3.6 of \cite{dette_optimal_2016-1}. The specifications of this four-parameter model are
\begin{eqnarray*}
	f^\T(x) & = & \left(1, x, x^2, x^3 \right), \quad x  \in  [1,2], \\
	\Cov\{\varepsilon(x),\varepsilon( x')\} & = & \min\left(x,x'\right), \qquad
	\lambda_{\min}(C)  =  0.0025.
\end{eqnarray*}

We again consider the D-criterion, so for $n = 5$, Table~\ref{tab:example2} contains the selected design points and D-efficiencies. This is one of the rare cases in which our method does even slightly better than the algorithm using the approximate sensitivity function from \cite{fedorov_design_1996}. All methods provide considerable improvements over uniform random designs. 

% latex table generated in R 4.0.2 by xtable 1.8-4 package
% Mon Oct 05 16:21:37 2020
\begin{table}[ht]
	\begin{center}
	\caption{Optimal designs and D-efficiencies for Example~2}
	\begin{tabular}{lcccccc}
			& $x_1$ & $x_2$ & $x_3$ & $x_4$ & $x_5$ & D-eff \\ 
			Q-VN & 1.00 & 1.16 & 1.52 & 1.84 & 2.00 & 0.9251 \\ 
			Q-VN+EP & 1.00 & 1.20 & 1.52 & 1.82 & 2.00 & 0.9300 \\ 
			Q-DPZ+EP & 1.00 & 1.14 & 1.33 & 1.60 & 2.00 & 0.8554 \\ 
			R-UNIF &  &  &  &  &  & 0.3208 \\ 
			BKSF & 1.00 & 1.16 & 1.46 & 1.83 & 2.00 & 0.9270 \\ 
			EXS & 1.00 & 1.21 & 1.61 & 1.84 & 2.00 & 0.9308 \\  
	\end{tabular}
	\end{center}
	\label{tab:example2}
\end{table}

\begin{figure}[hbtp!]
	\vskip-10pt
	\centering
	\begin{tabular}{cc}
		\includegraphics[width=0.45\textwidth]{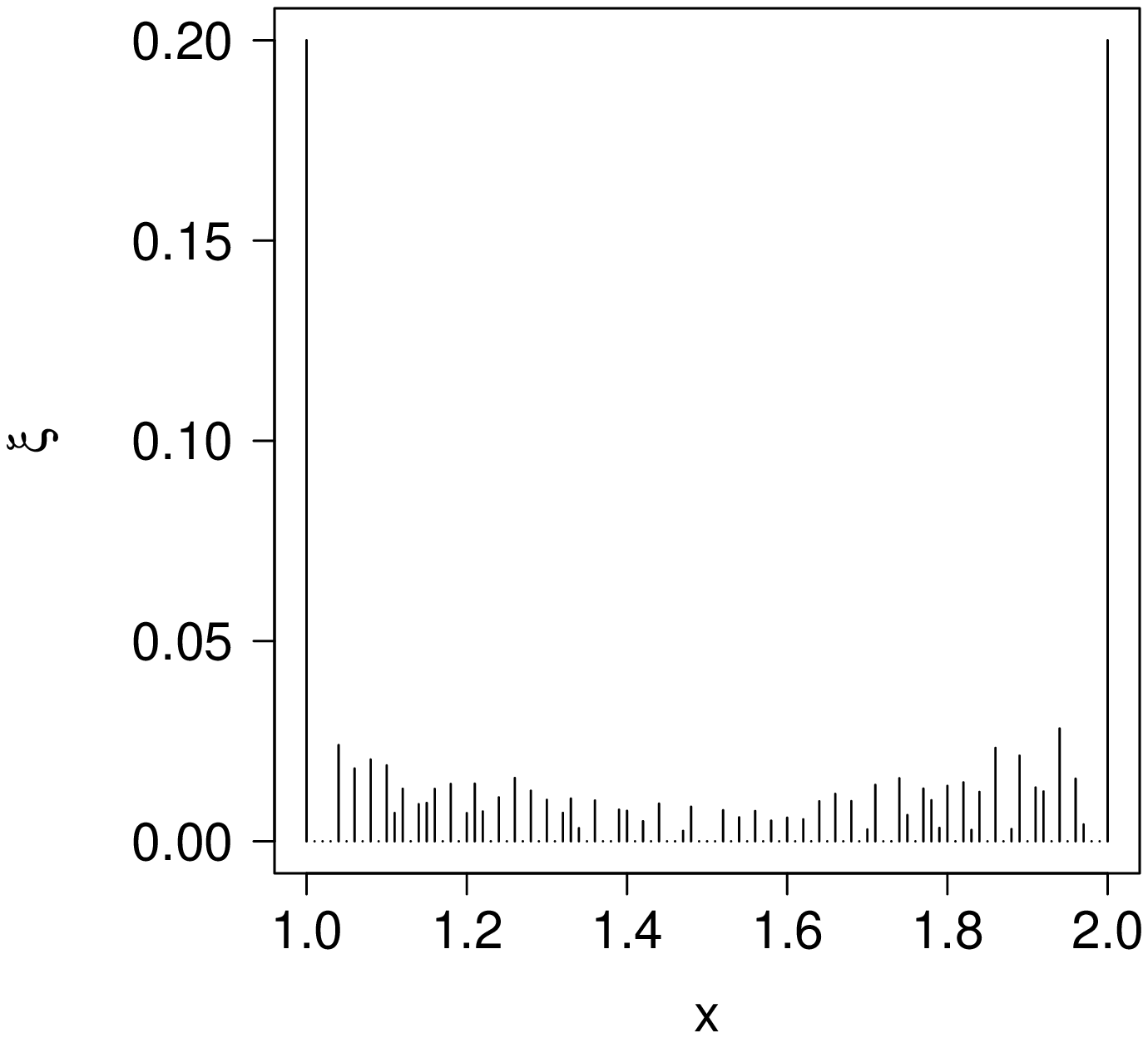} & 
		\includegraphics[width=0.45\textwidth]{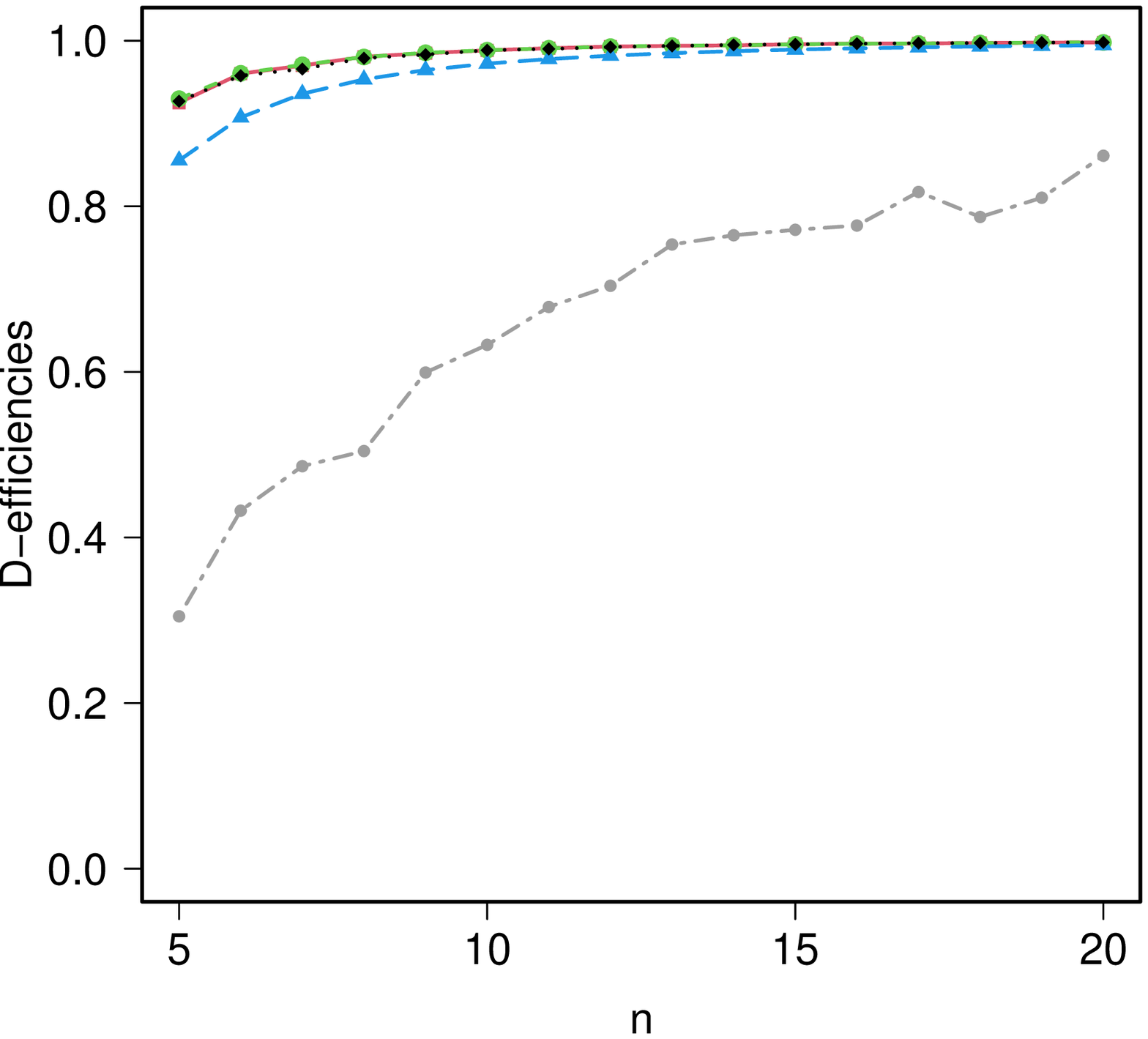}
	\end{tabular}
	\vskip-10pt
	\caption{Our measure (left panel) and efficiencies versus sample size (right panel) for Example~2. \label{fig:measures_example2}}
\end{figure}

%%%%%%%%%%%%%%%%%%%%%%%%%%%%%%%%%%%%%%%%%%%%%%%%%%%%
%%%%%%%%%%%%%%%%%%%%%%%%%%%%%%%%%%%%%%%%%%%%%%%%%%%%%%%%%%%%%%%%%%%%%%%%%%%%%%%%%%%%%%%%%%%%%

	\subsection{Example 3: a real application}\label{sec:example3}

	We consider a real-world application inspired by the example used throughout \cite{mateu_spatio-temporal_2012-1}. There a dataset with temperature and rainfall measurements gathered from 37 weather stations in the Austrian state of Upper Austria during the years 1994 -- 2009 was considered. Given this information, the goal was to find optimal designs for adding or reorganizing stations to optimize various aspects of future data collection.

	As for Figure 1.4 of \cite{mateu_spatio-temporal_2012-1} we employ the rainfall data from July 1994 and obtain the kriging estimates for an exponential kernel. The response function is assumed to be a plane, that is, it is linear in the coordinates. For our example, we assume that the parameters of the kernel function are known and set them to the kriging estimates. The model we consider is therefore
	\begin{eqnarray*}
		f^\T(x) & = & \left(1, x_1, x_2\right), \qquad x = (x_1,x_2)^\T, \\
		\Cov\{\varepsilon(x),\varepsilon( x')\} & = & 1756.65 \cdot \exp\left( \left\| x - x' \right\|_2 \bigl/ 40792.35 \right). 
	\end{eqnarray*}

	The design set is composed of the centroids of the $N = 445$ Upper Austrian municipalities, for which $\lambda_{\min}(C)  =  42.15$. We set $n = 36$ because this was the number of active weather stations in July 1994. Our hypothetical objective is therefore to reorganize the whole weather station network in order to most efficiently estimate the response function parameters given the assumed kernel function. We use the D-criterion as in the previous examples.

	The left panel of Figure~\ref{fig:measures_example3} displays the optimal measure obtained for the virtual noise representation when $n = 36$. The majority of the mass is concentrated on the borders. The right panel shows the optimal 36-point design we obtain by selecting the design with the highest criterion value out of 100 random designs sampled according to the measure depicted in the left panel. 

	\begin{figure}[hbtp!]
		\centering
		\begin{tabular}{cc}
			\includegraphics[width=0.45\textwidth]{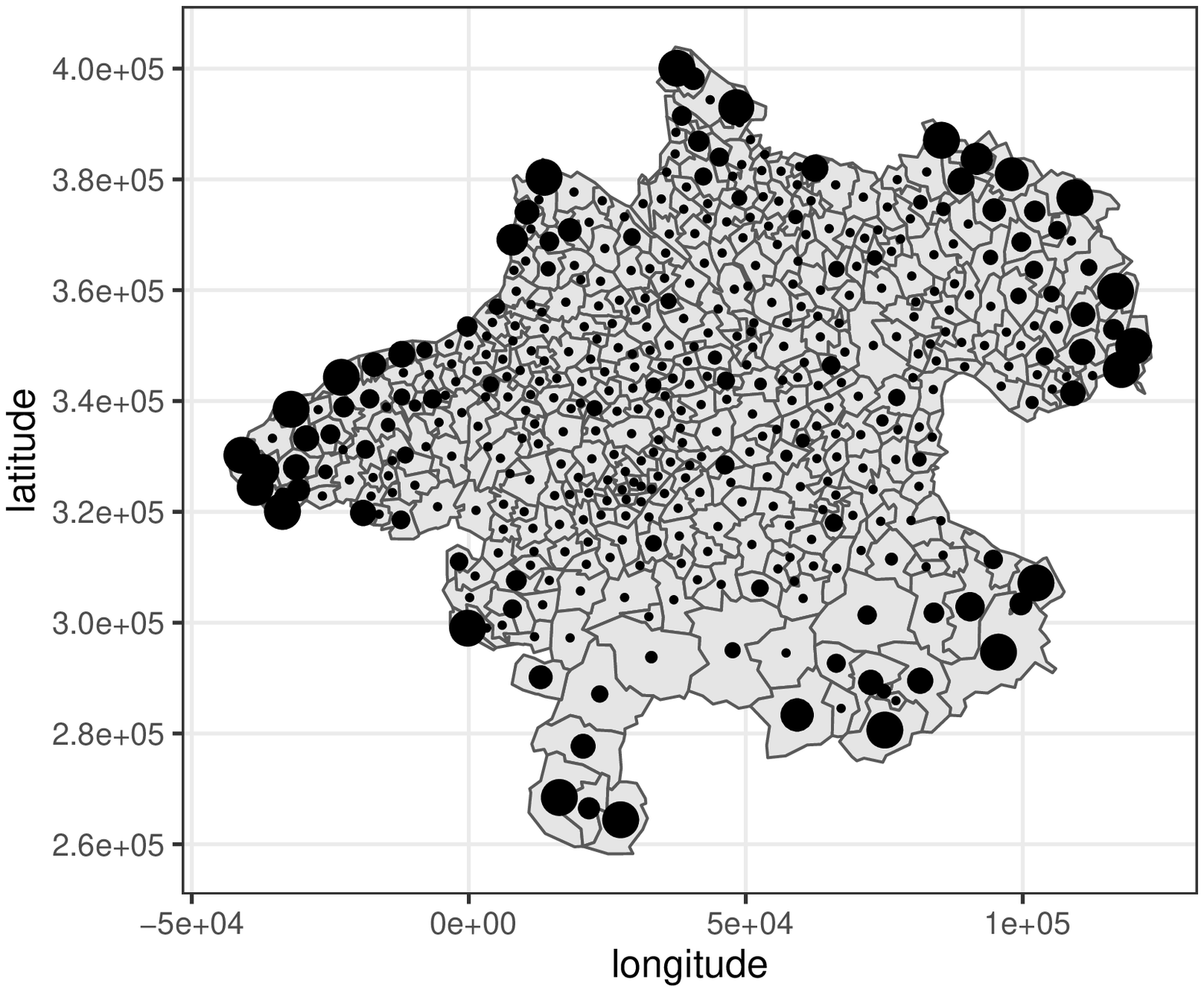} & 
			\includegraphics[width=0.45\textwidth]{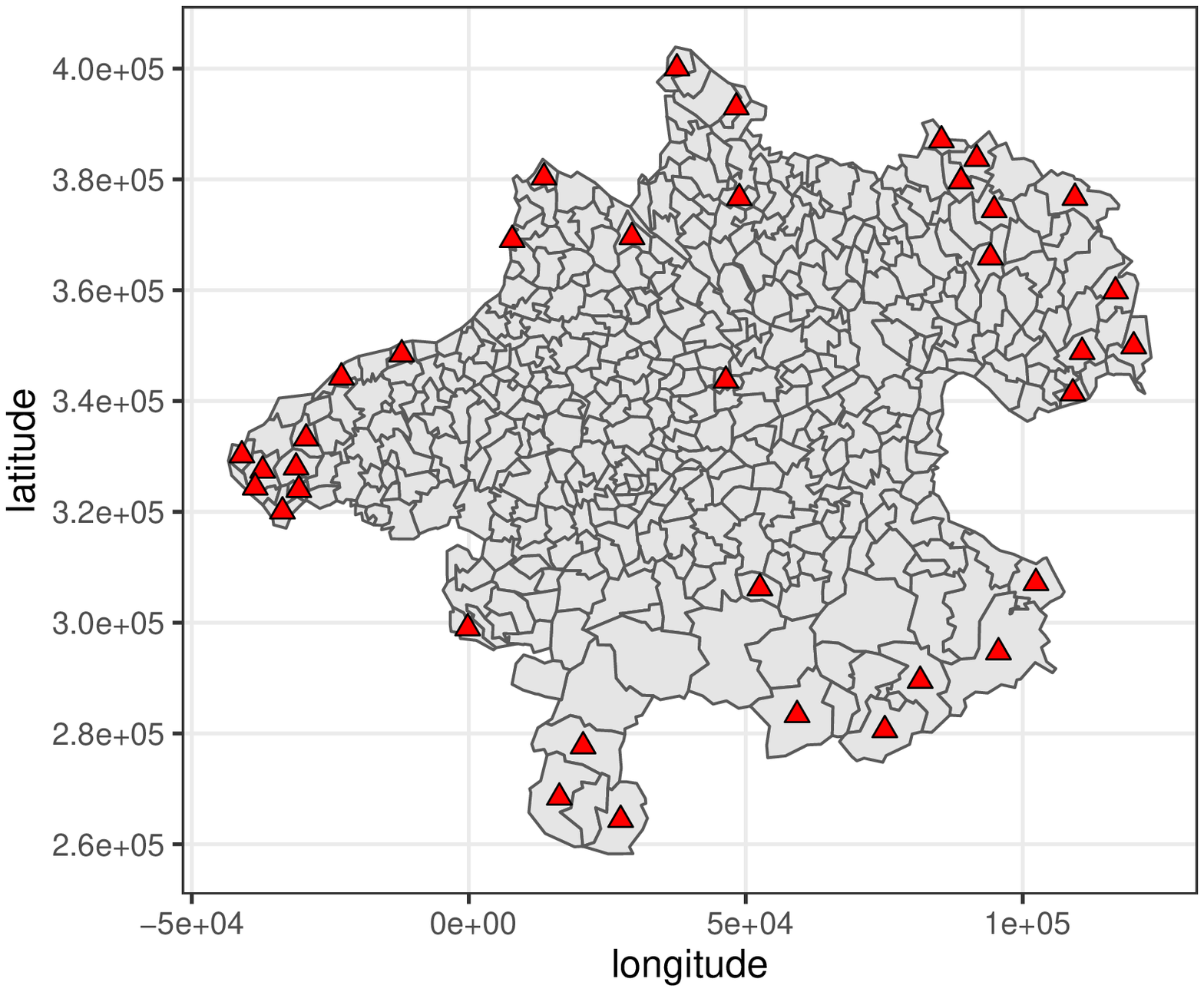}
		\end{tabular}
		\vskip-10pt
		\caption{Our measure (left panel) and exact 36-point design found by taking the optimal design from 100 random draws from the virtual noise measure (right panel) for Example~3 \label{fig:measures_example3}.}
	\end{figure}

	The D-efficiencies for different methods used to construct the optimal 36-point exact design are given in Table~\ref{tab:example3}. As in the other examples, these efficiencies are reported with respect to 
	our bound. The appendix contains those for additional values of $n$ from $n=4$ to $n=40$.

	\begin{table}[ht]
		\begin{center}
		\caption{D-efficiencies for Example~3}
		\begin{tabular}{lc}
				& D-eff \\ 
				R-VN (highest efficiency) & 0.9915 \\ 
				R-VN (median efficiency) & 0.9702 \\ 
				R-UNIF (highest efficiency) & 0.8405 \\ 
				R-UNIF (median efficiency) & 0.6689 \\ 
				BKSF & 0.9965 \\  
		\end{tabular}
		\label{tab:example3}
		\end{center}
	\end{table}
%%%%%%%%%%%%%%%%%%%%%%%%%%%%%%%%%%%%%%%%%%%%%%%%%%%%%%%%%%%%%%%%%%%%%%%%%%%%%%%%%%%%%%%%%%%%%

\section{Discussion}

The present paper complements and in a certain sense completes the research in \cite{muller_measures_2003}. The virtual noise approach now exhibits the convexity property it previously lacked. The formulated equivalence theorem allows for calculating a general upper bound for designs in Gaussian process regression, a key methodology in many application fields. This upper bound for the first time offers the possibility to calibrate and scale other design methods proposed in the literature.

Furthermore, our ``importance measures'' can be used to directly produce exact designs, be it by taking quantiles or by randomly sampling from them. However, typically for all the examples considered, also those not reported here, the adaptation of \cite{brimkulov_numerical_1980} using the approximate sensitivity function from \cite{fedorov_design_1996} performed slightly better than taking the quantiles from the measure given by the virtual noise representation, which in turn was better than the method suggested by \cite{dette_optimal_2016-1}.
Therefore, for practical purposes a combined approach using the approximate sensitivity function for finding designs and the upper bound derived from the virtual noise representation for evaluating them seems to be the most suitable avenue. Even when one prefers not to use an optimal design, but a classical one such as fractional factorial, orthogonal array, latin hypercube, etc. our bound can be used to judge how much efficiency would be sacrificed.

\section*{Acknowledgements} 
A. P\'azman was supported by the Slovak VEGA grant No. 1/0341/19.
M. Hainy was supported by the Austrian Science Fund (FWF): J3959-N32.
W.G. M{\"u}ller was partially supported by project grants LIT-2017-4-SEE-001 funded by the Upper Austrian Government, and Austrian Science Fund (FWF): I 3903-N32.

\renewcommand{\baselinestretch}{1}
\normalsize
\bibliographystyle{asa}
\bibliography{20pazman}

\newpage

\appendix

\section*{Appendix} 

%%%%%%  IMPORTANT %%%%%%%%%
\renewcommand{\thesection}{A.\arabic{section}}
% for arXive version use A., for paper supplement use S.

\section{Proof of Theorem 1
	on concavity}

	For any $u \in \mathbb{R}^N$ define the function
	$$
	\xi \in \Xi \to \gamma_u(\xi) = \left\{
	\begin{array}{ll}
	u^\T\left\{ C+W(\xi) \right\} ^{-1}u & \text{if } \supp(\xi)=\X \\
	u'^\T\left\{ C'+W'(\xi) \right\} ^{-1}u' &  \, \text{if } \supp(\xi)\neq \X  \\
	\end{array}
	\right. , 
	$$
	where $C'$, $W'(\xi)$, $u'$ are submatrices of $C$, $W(\xi)$, $u$, respectively with rows and/or columns corresponding to the points of $\supp(\xi)$.

\subsection*{Lemma A1}
	The function $\xi \in \Xi \to \gamma_u(\xi)$ is concave.

	\textit{Proof}.
	We start by proving the concavity on $\Xi_+=\left\{ \xi \in \Xi: \supp(\xi) = \X \right\}$. Consider an auxiliary function
	$$
	\psi_u ( \alpha ) :\alpha \in (0,1)\rightarrow \gamma_u(\xi_\alpha),
	$$
	with 
	$\xi _\alpha ( x ) =( 1-\alpha ) \xi ( x )
	+\alpha \mu ( x ); \quad \xi , \mu \in \Xi_+,
	\alpha \in (0,1)$. 
	The concavity of $\gamma_u(\xi)$ on $\Xi_+$ will be proven if we show that
 
\[
\frac{d^2 \psi_u ( \alpha ) }{d\alpha ^2} \le 0 
\]
for any $\xi, \mu \in \Xi_+ $ and any $\alpha \in \left( 0,1\right)$.

We can write 
\[
\psi_u ( \alpha ) =u^\T\left\{ H(\xi_\alpha) \right\} ^{-1}u,
\]
with 
\[
H(\xi_\alpha) =\left[ \left( C-\kappa I\right) +\frac \kappa n \, \diag\left\{ \frac 1{\xi _\alpha (
	\cdot ) }\right\}  \right].
\]
Notice that the matrix $\left( C-\kappa I\right)$ is positive definite,
hence also 
$H(\xi_\alpha)$ is positive definite.  

From the rule for the derivative of an inverse matrix we obtain
\begin{eqnarray*}
	\frac{d \psi_u ( \alpha ) }{d\alpha } &=&-\frac \kappa n \, u^\T H^{-1}(\xi_\alpha)  \left\{ \frac d{d\alpha } \, \diag\left( \frac 1{\xi
		_\alpha ( \cdot ) }\right) \right\} \, H^{-1}(\xi_\alpha) u
	\\
	&=&\frac \kappa n \, u^\T H^{-1}(\xi_\alpha) \, \diag\left\{ \frac
	1{\xi _\alpha ^2 ( \cdot ) }\right\}   \diag\left\{ \mu
	( \cdot ) -\xi ( \cdot ) \right\}  \, H^{-1}(\xi_\alpha) u \\
	&=&\frac \kappa n \, u^\T T^\T(\xi_\alpha) \, \diag\left\{ \mu
	( \cdot ) -\xi ( \cdot ) \right\}  \, T(\xi_\alpha) u,
\end{eqnarray*}
where 
\[
T(\xi_\alpha) =\left[   \left( C-\kappa I\right)\diag\left\{ \xi _\alpha (
\cdot ) \right\} +\frac \kappa n I\right]
^{-1}.
\]
Hence 
\begin{eqnarray*}
	&&\frac{d^2 \psi_u ( \alpha ) }{d\alpha ^2} \\
	&=&-\frac \kappa n \, u^\T T^\T(\xi_\alpha) \, \diag \{ \mu (
	\cdot ) -\xi ( \cdot ) \} \,\left( C-\kappa I\right)
	T^\T(\xi_\alpha) \\
	&&\times \, \diag \{ \mu (\cdot ) -\xi ( \cdot ) \} \,T(\xi_\alpha) u \\
	&&-\frac \kappa n \, u^\T T^\T(\xi_\alpha) \,\diag \{ \mu (\cdot ) -\xi ( \cdot ) \}  \\
	&&\times T(\xi_\alpha) \left( C-\kappa I\right) \diag \{ \mu (\cdot ) -\xi ( \cdot ) \} 
	\, T(\xi_\alpha) u,
\end{eqnarray*}
which is nonpositive since the matrix 
\begin{eqnarray*}
	\left( C-\kappa I\right) T^\T(\xi_\alpha) &=&\left( C-\kappa I\right)
	\left[ \diag\left\{ \xi _\alpha ( \cdot ) \right\}  \left( C-\kappa I\right) +\frac \kappa n I\right]^{-1}  \nonumber \\
	&=&\left[ \diag\left\{ \xi _\alpha ( \cdot ) \right\}
	+\frac \kappa n\left( C-\kappa I\right) ^{-1}\right] ^{-1}
\end{eqnarray*}
is symmetric positive definite.

To finish the proof take arbitrary $\xi,\mu \in \Xi$ and sequences $\{\xi_k\},\{\mu_k\}$ of elements of $\Xi_+$ converging to $\xi$ and $\mu$.
	From the inequality
	$$
	\gamma_u\{(1-\alpha)\xi_k + \alpha\mu_k\} \ge (1-\alpha)\gamma_u(\xi_k) + \alpha \gamma_u(\mu_k),
	$$
	valid for each $k$, we obtain by the continuity of $\xi \in \Xi \to \gamma_u(\xi)$, see Lemma A2, that the same inequality also holds for $\xi$ and $\mu$. 
\qed

\subsection*{Proof of Theorem 1
}

Since for any $t\in \mathbb{R}^p.$ we have 
	$$
	t^\T M(\xi) t = \gamma_{Ft}(\xi),
	$$
	we can write for any $\xi, \mu \in \Xi$
	$$
	t^\T M \{( 1-\alpha) \xi +\alpha \mu \} t \ge ( 1-\alpha )  t^\T M( \xi )t +\alpha  t^\T M( \mu) t,
	$$
which means that
\[
M\{ ( 1-\alpha ) \xi +\alpha \mu \} \geq ( 1-\alpha
) M ( \xi ) +\alpha M( \mu ) 
\]
in the L\"owner ordering. From the monotonicity of $\Phi $ we have 
\[
\Phi \left[ M\{ ( 1-\alpha ) \xi +\alpha \mu \}
\right] \geq \Phi \left\{ ( 1-\alpha ) M( \xi )
+\alpha M( \mu ) \right\} 
\]
and from the concavity of $\Phi \,$ we have 
\[
\Phi \left\{ ( 1-\alpha ) M( \xi ) +\alpha M( \mu ) \right\}  \geq ( 1-\alpha ) \Phi \{ M( \xi
) \} +\alpha  \Phi \{ M( \mu ) \}.
\]
\qed

\section{Proof of the equivalence theorem (Theorem~2)}

Since the function 
\[
\xi \in \Xi \rightarrow \Phi \{ M ( \xi ) \} 
\]
is concave, the design $\bar{\xi}$ is $\Phi -$optimal if and only if for
every $\mu \in \Xi $ we have, see Lemma A3, 
\begin{equation*}\label{A1}
\pi\left(\bar{\xi}, \mu\right)  =  \lim_{\alpha \rightarrow 0^+} \frac{\Phi \left\{ M\left( \xi_{\alpha}^{\mu} \right) \right\} - \Phi \left\{ M\left( \bar{\xi} \right) \right\} }{
	\alpha } \: \leq \: 0,
\end{equation*}
where we used the notation 
\[
\xi _\alpha ^\mu ( x ) =( 1-\alpha ) \bar{\xi}(
x) +\alpha \mu ( x ) .
\]

We have
\begin{eqnarray*}
\lim_{\alpha \rightarrow 0^+} \frac{\Phi \left\{ M\left( \xi_{\alpha}^{\mu} \right) \right\} - \Phi \left\{ M\left( \bar{\xi} \right) \right\} }{ \alpha } & = & 
\sum_{i,j=1}^p \lim_{\alpha \rightarrow 0^+} \frac{\Phi \left\{ M\left( \xi_{\alpha}^{\mu} \right) \right\} - \Phi \left\{ M\left( \bar{\xi} \right) \right\} }{ M_{ij} \left( \xi_{\alpha}^{\mu} \right) - M_{ij} \left( \bar{\xi} \right)} \times \\
& & \phantom{\sum_{i,j=1}^p} \times \lim_{\alpha \rightarrow 0^+} \frac{ M_{ij} \left( \xi_{\alpha}^{\mu} \right) - M_{ij} \left( \bar{\xi} \right)}{\alpha} \\
& = & 
\sum_{i,j=1}^p \lim_{\xi \rightarrow \bar{\xi}} \, \frac{\Phi \left\{ M\left( \xi \right) \right\} - \Phi \left\{ M\left( \bar{\xi} \right) \right\} }{ M_{ij} \left( \xi \right) - M_{ij} \left( \bar{\xi} \right)} \times \\
& & \phantom{\sum_{i,j=1}^p} \times \lim_{\alpha \rightarrow 0^+} \frac{ M_{ij} \left( \xi_{\alpha}^{\mu} \right) - M_{ij} \left( \bar{\xi} \right)}{\alpha},
\end{eqnarray*}
since $\underset{\alpha \rightarrow 0^+}{\lim} \, \xi_{\alpha}^{\mu} = \bar{\xi}$. Therefore, we obtain
\begin{eqnarray*}
	\lim_{\alpha \rightarrow 0^+} \frac{\Phi \left\{ M\left( \xi_{\alpha}^{\mu} \right) \right\} - \Phi \left\{ M\left( \bar{\xi} \right) \right\} }{ \alpha } & = & 
	\tr\left[ \nabla \Phi \left(M\right)\bigl|_{M = M\left(\bar{\xi}\right)} \: \lim_{\alpha \rightarrow 0^+} \frac{ M \left( \xi_{\alpha}^{\mu} \right) - M \left( \bar{\xi} \right)}{\alpha} \right] \\
	& = & 
	\tr\left[ \nabla \Phi \left(M\right)\bigl|_{M = M\left(\bar{\xi}\right)} \: \frac{\kappa}{n} F^{\T} T^{\T} \left( \bar{\xi }\right) \diag\left( \mu - \bar{
		\xi}\right)  T\left(\bar{\xi}\right) F \right] \\
	& = & \frac \kappa n\sum_{x\in \mathcal{X}}\left\{ \mu
	( x ) -\bar{\xi} ( x ) \right\} h(x,\bar\xi),
\end{eqnarray*}
as follows from \eqref{eq:dL_dalpha} and \eqref{eq:lim_M_alpha0} in Lemma~A4 and from the definition of $h(x,\bar{\xi})$ at the beginning of Section~\ref{subsec:equi_theorem}, just before Theorem~2. Hence we finally obtain
\begin{equation} \label{eq:pi}
	\pi\left(\bar{\xi},\mu\right) = \frac \kappa n\sum_{x\in \mathcal{X}}\left\{ \mu
	( x ) -\bar{\xi} ( x ) \right\} h(x,\bar\xi).
\end{equation}
Therefore, $\bar{\xi} \in \Xi$ is $\Phi -$optimal if and only if for every $\mu \in \Xi $
we have 
\begin{eqnarray} \label{A5}
&&\sum_{x\in \mathcal{X}}\mu ( x ) h(x,\bar\xi) 
\leq \sum_{x\in \mathcal{X}}\bar{\xi} ( x ) h(x,\bar\xi).
\end{eqnarray}
Since this expression is linear in $\mu $ and since the convex hull of
the set of all exact $n$-point
designs is the convex set 
$\Xi $, we do not need to consider all designs $\mu $ from $\Xi,$ but
just the exact $n$-point designs. Hence the validity of the inequality \eqref{equiv1} in Theorem~2 
for all such
$n$-point measures is equivalent to the validity of (\ref{A5}) for all $\mu \in \Xi $.

From \eqref{eq:pi} it follows that the inequalities \eqref{equiv1_delta} in Theorem~2 imply that for every $\mu \in \Xi$ we have
\[ \pi\left(\bar{\xi}, \mu\right) \leq \delta. \]
According to Lemma~A3 it follows that the inequality \eqref{equiv1_phi} holds as well. \qed

\section{Adaptation of the algorithm from \cite{fedorov_design_1996} for A-optimality}

\subsection*{Sensitivity function for the D-criterion}

Consider the linear regression model with i.i.d. errors,
\begin{equation*}
y(x) = f(x)^\T \theta + \varepsilon(x),
\end{equation*}
where $E\{\varepsilon(x)\} = 0$, $\Var\{\varepsilon(x)\} = \sigma^2_\varepsilon$, and $E\{\varepsilon(x) \varepsilon(x')\} = 0$ $\forall x, x' \in \mathcal{X}$.

The information matrix for a design $\xi \in \Xi$ is
\begin{equation*}
M_{\xi} = \frac{1}{\sigma^2_\varepsilon} \sum_{x \in \mathcal{X}} \xi(x) f(x) f(x)^\T.
\end{equation*}

We also consider a one-point design at $\bar{x} \in \mathcal{X}$ denoted by $\xi_{\bar{x}}$. Its information matrix is
\begin{equation*}
M_{\xi_{\bar{x}}} = \frac{1}{\sigma^2_\varepsilon} f(\bar{x}) f(\bar{x})^\T.
\end{equation*}

Let the information matrix of the convex combination $\xi_{\bar{x}}^{\alpha} = (1 - \alpha) \xi + \alpha \xi_{\bar{x}}$ of those two designs be denoted by
\begin{equation*}
M_{\alpha} = (1 - \alpha) M_{\xi} + \alpha M_{\xi_{\bar{x}}}.
\end{equation*}

For some concave criterion function $\Phi(M)$, the directional derivative of $\Phi(\cdot)$ at the design $\xi$ in the direction of $\xi_{\bar{x}}$, the so-called sensitivity function, is 
\begin{eqnarray}
\phi(\bar{x}; \xi) & = & \lim_{\alpha \rightarrow 0^+} \frac{\Phi(M_{\alpha}) - \Phi(M_{\xi})}{\alpha} \notag \\
& = & \lim_{\alpha \rightarrow 0^+} \tr \left\{ \nabla_M \Phi(M_{\alpha}) \, \frac{\partial M_{\alpha}}{\partial \alpha} \right\} \notag \\
& = &  \lim_{\alpha \rightarrow 0^+} \tr \left\{ \nabla_M \Phi(M_{\alpha}) \left( M_{\xi_{\bar{x}}} - M_{\xi} \right) \right\} \notag\\
& = & \tr \left\{ \nabla_M \Phi(M_{\xi}) \left( M_{\xi_{\bar{x}}} - M_{\xi} \right) \right\} \notag \\
& = & \frac{1}{\sigma^2_\varepsilon} f(\bar{x})^\T \nabla_M \Phi(M_{\xi}) f(\bar{x}) - \tr \left\{ \nabla_M \Phi(M_{\xi}) \, M_{\xi} \right\}. \label{eq:sensitivity_function_general}
\end{eqnarray}

In the case of \textbf{D-optimality}, we prefer here the criterion function $\Phi_D(M) = \log \det (M)$, so $\nabla_M \Phi_D(M) = M^{-1}$. Plugging this into \eqref{eq:sensitivity_function_general} yields
\begin{eqnarray}
\phi_D(\bar{x}; \xi) & = & \frac{1}{\sigma^2_\varepsilon} f(\bar{x})^\T M_{\xi}^{-1} f(\bar{x}) - \tr \left(M_{\xi}^{-1} M_{\xi} \right) \notag \\
& = & \frac{1}{\sigma^2_\varepsilon} f(\bar{x})^\T M_{\xi}^{-1} f(\bar{x}) - \mathrm{dim}(\theta). \label{eq:sensitivity_function_D}
\end{eqnarray}

Up to a constant, expression \eqref{eq:sensitivity_function_D} is equal to 
\begin{equation}
\frac{E\left[ \left\{y(\bar{x}) - \hat{y}(\bar{x})\right\}^2 \right]}{\sigma^2_\varepsilon} = \frac{\sigma^2_\varepsilon + f(\bar{x})^\T M_{\xi}^{-1} f(\bar{x})}{\sigma^2_\varepsilon} = 1 + \frac{f(\bar{x})^\T M_{\xi}^{-1} f(\bar{x})}{\sigma^2_\varepsilon}. \label{eq:Brimkulov_function_D}
\end{equation}

In the correlated case with no replications, \cite{fedorov_design_1996} suggests to replace the numerator in (\ref{eq:Brimkulov_function_D}) by the unconditional and the denominator in (\ref{eq:Brimkulov_function_D}) by the conditional variance. 
Let $\Tc = \{x_1,\ldots,x_n\} \subset \X$ denote an exact n-point design, let $C(\Tc)$ be the $n \times n$ submatrix of the covariance matrix $C$ for the set of points in $\Tc$, and let $F(\Tc)$ be the $n \times p$ design matrix with elements $F_{i,j}(\Tc) = f_j(x_i)$. Let $k(x,x')=\Cov\{\varepsilon(x),\varepsilon( x')\}$ and for a point $\bar{x} \in \X \backslash \Tc$, let $k(\bar{x},\Tc) = \left(k(\bar{x},x_1),\ldots,k(\bar{x},x_n) \right)^\T$. In the correlated case, \cite{fedorov_design_1996} therefore replaces expression~(\ref{eq:Brimkulov_function_D}) with

\begin{equation*}
\frac{E\left[\left\{ y(\bar{x}) - \hat{y}(\bar{x})\right\}^2 \right]}{\tilde{\sigma}^2(\bar{x})} = \frac{\tilde{\sigma}^2(\bar{x}) + \tilde{f}(\bar{x})^\T M_{\Tc}^{-1} \tilde{f}(\bar{x})}{\tilde{\sigma}^2(\bar{x})},
\end{equation*}
where
	\begin{equation*}
	\tilde{\sigma}^2(\bar{x}) = k(\bar{x},\bar{x}) - k(\bar{x},\Tc)^\T C^{-1}(\Tc) k(\bar{x},\Tc), \label{eq:sigma2_tilde}
	\end{equation*}
	is the conditional variance,  $\tilde{f}(\bar{x})$ is
\begin{equation*}
\tilde{f}(\bar{x}) = f(\bar{x}) - F(\Tc)^\T C^{-1}(\Tc) k(\bar{x},\Tc), \label{eq:f_tilde}
\end{equation*}
and the information matrix is
\begin{equation*}
M_{\Tc} = F(\Tc)^\T C^{-1}(\Tc) F(\Tc). \label{eq:M_T}
\end{equation*}

It can be shown that for D-optimality this approximate sensitivity function is equal to the factor by which the determinant of the information matrix is increased by adding design point $\bar{x}$.

\subsection*{Sensitivity function for the A-criterion}

In the case of \textbf{A-optimality}, we use the criterion $\Phi_A(M) = -\mathrm{tr}\left(M^{-1}\right)$. We therefore have $\nabla_M \Phi_A(M) = -\left(-M^{-2}\right)^\T = M^{-2}$ and the directional derivative is
\begin{equation*}
\phi_A(\bar{x}; \xi) = \frac{1}{\sigma^2_\varepsilon} f(\bar{x})^\T M_{\xi}^{-2} f(\bar{x}) - \mathrm{tr}\left( M_{\xi}^{-1}\right). \label{eq:sensitivity_function_A}
\end{equation*}

In accordance with the case of D-optimality outlined above, we obtain the approximate sensitivity function by analogous replacements.
Therefore, the sensitivity function we eventually use in some of our examples is
\begin{equation*}
\tilde{\phi}_A(\bar{x}; \Tc) = \frac{1}{\tilde{\sigma}^2(\bar{x})} \tilde{f}(\bar{x})^\T M_{\Tc}^{-2} \tilde{f}(\bar{x}) - \mathrm{tr}\left(M_{\Tc}^{-1}\right).
\end{equation*}

\section{Outline of algorithm from \cite{fedorov_design_1996}}

\newcommand{\designgrid}{\mathcal{X}}
\newcommand{\candidateset}{\mathcal{C}}

Algorithm~\ref{algo:BKSF} contains the description of our implementation for a finite design grid of the algorithm proposed in \cite{fedorov_design_1996}.

\begin{algorithm}[H]
	\caption{Algorithm proposed by \cite{fedorov_design_1996} for a finite design grid\label{algo:BKSF}}
	\SetKw{TRUE}{true}
	\SetKw{FALSE}{false}
	\SetKw{NOT}{not}
	\KwIn{Design grid $\designgrid$ of possible design points; design size $n$; initial design $\Tc_0$; (approximate) sensitivity function $\tilde{\phi}(x;\Tc).$}
	\KwOut{Exact design $\Tc$ found by algorithm.}
	
	\texttt{abort} = \FALSE\;
	$\Tc = \Tc_0$\;
	\While{\NOT \texttt{abort}}{ 
		\For{$i = 1$ \KwTo $n$}{ 
			Let $x_i$ be the $i$\textsuperscript{th} design point in $\Tc$\;
			$\Tc_{-i} = \Tc \backslash \{x_i\}$\;
			$\mathtt{s[i]} = \tilde{\phi}(x_i;\Tc_{-i})$\;
		}
		$\mathtt{sens\_drop} = \min(\mathtt{s})$\;
		Let $\mathtt{k}$ be the index where $\mathtt{s[k]} = \mathtt{sens\_drop}$\;
		Let $\candidateset = \designgrid \backslash \Tc_{-k}$\;
		$m = \text{card}(\candidateset)$\;
		\For{$i = 1$ \KwTo $m$}{
			Let $\bar{x}_i$ be the $i$\textsuperscript{th} design point in the candidate set $\candidateset$\;
			$\mathtt{t[i]} = \tilde{\phi}(\bar{x}_i; \Tc_{-k})$\;
		}
		$\mathtt{sens\_add} = \max(\mathtt{t})$\;
		Let $\mathtt{l}$ be the index where $\mathtt{t[l]} = \mathtt{sens\_add}$\;
		$\Tc = \Tc_{-k} \cup \bar{x}_l$\;
		\If{$\mathtt{sens\_add} - \mathtt{sens\_drop} \leq 0$}{
			\texttt{abort} = \TRUE\;
		}
	}
\end{algorithm}

\section{Miscellanea}

	\subsection*{Lemma A2: the continuity lemma}
	The function $\gamma_u(\xi)$ defined in Lemma A1 is continuous on $\Xi$ and $M(\xi)$ is continuous on $\Xi$ as well.
	
	\textit{Proof}.
	Take $\xi_k \in \Xi$ such that $\lim_{k \to \infty} \xi_k = \xi$. The equality $\lim_{k \to \infty} \gamma_u(\xi_k) = \gamma_u(\xi)$ is evident if $\supp(\xi)=\X$, or more generally if $\bigcup_{k_0=1}^\infty \bigcap_{k \geq k_0} \supp(\xi_k)=\supp(\xi)$. To clarify how to do the proof in the opposite case we suppose that  $\supp(\xi)=\X$, but $\bigcup_{k_0=1}^\infty \bigcap_{k \geq k_0} \supp(\xi_k)$ has one point less than the set $\supp(\xi)$, say the point $x_1$. By a basic property of inverse matrices we have
	\begin{equation}\label{cont}
	\left[\left\{C+W(\xi_k)\right\}^{-1}\right]_{ij} = \frac{(-1)^{i+j}\det\left\{C+W(\xi_k)\right\}^{i,j}}{\det\left\{C+W(\xi_k)\right\}},
	\end{equation}
	where $\left\{C+W(\xi_k)\right\}^{i,j}$ is the submatrix of $\left\{C+W(\xi_k)\right\}$ after omitting the $i$-th row and $j$-th column. For $k\to\infty$ (\ref{cont}) tends to zero for $i=1$ or $j=1$, but to a nonzero number when $i\neq1$ and $j\neq1$. This can be seen when using the definition of a determinant as a sum of products of the elements of the matrix. When $k\to\infty$, some terms of these sums in the numerator and the denominator converge to $\infty$ much slower than the others, so we can neglect them.
	We proceed similarly when
	$\bigcup_{k_0=1}^\infty \bigcap_{k \geq k_0} \supp(\xi_k) \cup \{x_1,x_2\}=\supp(\xi)$, etc.. \qed

\vspace*{2ex}

For the reader's convenience we prove a statement from convex theory, see
e.g. \cite{rockafellar_convex_1970}.

\subsection*{Lemma A3}
If a function $\phi :\xi \in \Xi \rightarrow \phi(\xi) \in \mathbb{R}$ is concave, then for any 
$\bar{\xi},\mu \in \Xi $ the function 
\[
\alpha \in \left( 0,1\right) \rightarrow \frac{\phi \left\{ ( 1-\alpha
	) \bar{\xi}+\alpha \mu \right\} -\phi \! \left( \bar{\xi}\right) }\alpha 
\]
is nonincreasing. Hence for any $\bar{\xi} \in \Xi$ and any $\mu \in \Xi$ there exists the
limit
\[
\lim_{\alpha \rightarrow 0^{+}}\frac{\phi \left\{ ( 1-\alpha ) 
	\bar{\xi}+\alpha \mu \right\} -\phi \! \left( \bar{\xi}\right) }\alpha \equiv \pi\left(\bar{\xi}, \mu \right), 
\]
and we have
\[
\bar{\xi}\in \arg \max_{\xi \in \Xi }\phi ( \xi ) 
\]
if and only if 
\begin{equation} \label{A6}
\forall \: \mu \in \Xi \,\,\,\,\,\,\,\,\, \pi\left(\bar{\xi}, \mu \right) \: \leq \: 0. \quad \quad
\end{equation}

	If for some $\delta > 0$ and for every $\mu \in \Xi$ we have
	$$ 
	 \pi\left(\bar{\xi}, \mu \right) \leq \delta,
	$$
	then
	$\displaystyle \phi\left(\bar\xi\right) \ge \underset{\xi \in \Xi}{\max}\ \phi(\xi)-\delta$.

\subsection*{Proof}

Take $\,0<\alpha _1<\alpha _2<1$. From the concavity of $\phi
( \xi ) $ we obtain 
\begin{eqnarray*}
	&&\phi \left\{ ( 1-\alpha _1 ) \bar{\xi}+\alpha _1\mu \right\} -\phi \!
	\left( \bar{\xi}\right) \\
	&=&\phi \left[ \frac{\alpha _1}{\alpha _2}\left\{ ( 1-\alpha _2 ) 
	\bar{\xi}+\alpha _2\mu \right\} +\left( 1-\frac{\alpha _1}{\alpha _2}\right) 
	\bar{\xi}\right] -\phi \! \left( \bar{\xi}\right) \\
	&\geq &\frac{\alpha _1}{\alpha _2} \, \phi \left\{ ( 1-\alpha _2 ) \bar{
		\xi}+\alpha _2\mu \right\} +\left( 1-\frac{\alpha _1}{\alpha _2}\right) \phi \!
	\left( \bar{\xi}\right) -\phi \! \left( \bar{\xi}\right) \\
	&=&\frac{\alpha _1}{\alpha _2}\left[ \phi \left\{ ( 1-\alpha _2 ) 
	\bar{\xi}+\alpha _2\mu \right\} -\phi \! \left( \bar{\xi}\right) \right].
\end{eqnarray*}
We multiply this inequality by $1 / \alpha_1$ to prove the first
statement of the Lemma.

If $\bar{\xi}\in \underset{\xi \in \Xi}{\arg \max}\ \phi ( \xi ),$ then for
any $\mu $ and $\alpha $ we evidently have
\[
\frac{\phi \left\{ ( 1-\alpha ) \bar{\xi}+\alpha \mu \right\} -\phi \! 
	\left( \bar{\xi}\right) }\alpha \leq 0
\]
and this inequality does also hold for $\alpha \rightarrow 0^{+}.$

On the other hand, suppose that there is $\mu \in \Xi $ such that $\phi
(\mu) >\phi \! \left( \bar{\xi}\right) $. Then for every $\alpha
	\in \left( 0, 1\right) $ we have 
\begin{eqnarray*}
	&&\frac{\phi \left\{ ( 1-\alpha ) \bar{\xi}+\alpha \mu \right\}
		-\phi \! \left( \bar{\xi}\right) }\alpha  \\
	&\geq &\frac{( 1-\alpha ) \phi \! \left( \bar{\xi}\right) +\alpha
		\phi \! \left( \mu \right) -\phi \! \left( \bar{\xi}\right) }\alpha  \\
	&=&\phi \! \left( \mu \right) -\phi \! \left( \bar{\xi}\right) > 0,
\end{eqnarray*}
and by taking  the limit for $\alpha \rightarrow 0$ we see that (\ref{A6}) does not hold. 

	If 
	$$ 
	\forall_{\mu \in \X} \: \pi\left(\bar{\xi}, \mu \right) \leq \delta,
	$$
	then from the concavity of $\phi(\xi)$ we obtain for every $\mu \in \Xi$
	\begin{eqnarray*}
		\phi(\mu)-\phi(\bar\xi) &=& \frac{(1-\alpha)\phi(\bar\xi)+\alpha\phi(\mu)-\phi(\bar\xi)}{\alpha} 
		\\
		&\leq& 
		\frac{\phi\{(1-\alpha)\bar\xi+\alpha\mu\}-\phi(\bar\xi)}{\alpha} \leq   \pi\left(\bar{\xi}, \mu \right) \leq \delta.
	\end{eqnarray*}
	Hence $\displaystyle \phi\left(\bar\xi\right) \ge \underset{\mu \in \Xi}{\max}\ \phi(\mu)-\delta$.
\qed

\vspace*{2ex}

\subsection*{Lemma A4}

The matrix 
\begin{equation}
L\left( \xi \right) \equiv F^\T \left[ \diag\left\{\xi \left( \cdot \right) \right\}
\left(C-\kappa I\right) +\frac{\kappa}{n} I\right]^{-1} \diag\left\{ \xi \left( \cdot \right) \right\} F
\label{eq:L_xi}
\end{equation}
is well defined and continuous on $\Xi$.

For every $\xi \in \Xi $ we have $L\left( \xi \right) = M\left( \xi \right)$
with $M\left( \xi \right)$ defined in \eqref{FIMvirtual} and \eqref{FIMvirtual2}.

For every $\bar{\xi} \in \Xi$ and $\mu \in \Xi$ the right-hand
limit 
\[
\lim_{\alpha \rightarrow 0^{+}}\frac{L\left\{ \left( 1-\alpha \right) 
	\bar{\xi} +\alpha \mu \right\} -L\left( \bar{\xi} \right) }
{\alpha}
\]
is well defined and it is a continuous function of $\bar{\xi}$ on
the whole set $\Xi$. On the set $\Xi _{+}$ it is equal to the derivative $
\frac{\partial L\left( \xi _\alpha ^\mu \right) }{\partial \alpha }\bigl|
_{\alpha =0}=\frac{\partial M\left( \xi _\alpha ^\mu \right) }{\partial
	\alpha }\bigl|_{\alpha =0}$.

\subsection*{Proof}

We denote $Z\left( \xi \right) = \left[ \diag\left\{\xi \left( \cdot \right) \right\}
\left(C-\kappa I\right) +\frac{\kappa}{n} I\right]$. To see the
correctness of the definition and the continuity of $L\left( \xi \right)$,
we must prove the existence of the matrix $Z\left( \xi \right) ^{-1}$ used in 
\eqref{eq:L_xi}. Without loss of generality we can suppose that $\xi
\left( x_i\right) >0$ for $i=1,\ldots,k$ and $\xi \left( x_i\right) = 0$ for $
i=k+1,...,N$, where $k\leq N$. Let us write $Z\left( \xi \right) =\left( 
\begin{array}{cc}
V & R \\ 
P & T
\end{array}
\right)$, where $V$ is a $k\times k$ matrix and $T$ is a $\left( N-k\right)
\times \left( N-k\right) $ matrix. Evidently $P=0.$ According to \cite{harville_matrix_algebra_1997}, Lemma 8.5.4, if the inverse matrices $V^{-1}$ and $T^{-1}$ exist,
then the inverse of the partitioned matrix exists as well, and 
\[
\left( 
\begin{array}{cc}
V & R \\ 
0 & T
\end{array}
\right) ^{-1}=\left( 
\begin{array}{cc}
V^{-1} & -V^{-1}RT_{}^{-1} \\ 
0 & T_{}^{-1}
\end{array}
\right). 
\]

We can write 
\[
\diag\left\{ \xi \left( \cdot \right) \right\} =\left( 
\begin{array}{cc}
\diag' \left\{ \xi \left( \cdot \right) \right\} & 0 \\ 
0 & 0
\end{array}
\right),
\]
where the $k\times k$ matrix $\diag^{\prime }\left\{ \xi \left( \cdot \right)
\right\} $ has the $k$ positive elements of $\xi $ on its
diagonal. Consequently 
\begin{eqnarray*}
	V &=& \diag^{\prime}\left\{ \xi \left( \cdot \right) \right\} \left( C^{\prime
	}-\kappa I^{\prime}\right) +\frac \kappa nI^{\prime } \\
	&=& \diag^{\prime}\left\{ \xi \left(\cdot\right) \right\} \left[ \left(
	C^{\prime}-\kappa I^{\prime}\right) +\frac \kappa n \diag^{\prime
	}\left\{ \xi ^{-1}\left( \cdot\right) \right\} \right],
\end{eqnarray*}
where the prime denotes the $k\times k$ submatrix of the corresponding
matrix. The symmetric matrix
\[
\left[ \left(
C^{\prime}-\kappa I^{\prime}\right) +\frac \kappa n \diag^{\prime
}\left\{ \xi ^{-1}\left( \cdot\right) \right\} \right]
\]
is positive definite, hence its inverse exists and $V^{-1}$ exists as
well. Finally $T=\frac \kappa nI^{*}$, where the star denotes the $\left(
N-k\right) \times \left( N-k\right) $ submatrix, so $T^{-1}=\frac n\kappa
I^{*}.$

The continuity of $L\left( \xi \right) $ on $\Xi $ is now evident from \eqref{eq:L_xi}.

To prove $L\left( \xi \right) =M\left( \xi \right)$, let us denote the upper $k\times p$ and lower $\left(
N-k\right) \times p$ submatrices of $F$ by $F^{\prime}$ and $F^{*}$, respectively. We have 
\begin{eqnarray*}
	L\left( \xi \right) &=&\left( 
	\begin{array}{c}
		F^{\prime} \\ 
		F^{*}
	\end{array}
	\right)^\T \left( 
	\begin{array}{cc}
		V^{-1} & -V^{-1}RT_{}^{-1} \\ 
		0 & T_{}^{-1}
	\end{array}
	\right) \left( 
	\begin{array}{c}
		\diag^{\prime}\left\{ \xi \left( \cdot \right) \right\} F^{\prime} \\ 
		0
	\end{array}
	\right) \\
	&=&\left( F^{\prime}\right)^\T V^{-1} \diag^{\prime}\left\{ \xi \left(
	\cdot\right) \right\} F^{\prime} \\
	&=&\left( F^{\prime}\right)^\T \left[ \left(
	C^{\prime}-\kappa I^{\prime}\right) +\frac \kappa n \diag^{\prime
	}\left\{ \xi ^{-1}\left( \cdot\right) \right\} \right]^{-1} F^{\prime} = M\left( \xi
	\right),
\end{eqnarray*}
as follows from Eq.~\eqref{FIMvirtual2}.

Let us suppose first that $\bar{\xi} \in \Xi_{+}$. Using the
standard rules for the derivatives of matrices, we obtain from \eqref{eq:L_xi}
\begin{eqnarray}
	\frac{\partial L\left( \xi _\alpha ^\mu \right) }{\partial \alpha }\biggl|
	_{\alpha =0} & = & - \, F^{\T} Z^{-1}\left( \bar{\xi}\right) \, \diag\left( \mu 
	- \bar{\xi }\right) \left( C-\kappa I\right) Z^{-1}\left( \bar{\xi }\right) \diag\left( \bar{\xi } \right) F \notag \\
	&& + \, F^{\T} Z^{-1}\left( \bar{\xi }\right) \diag\left( \mu - \bar{
		\xi}\right) F \notag \\
&=& F^{\T} Z^{-1}\left( \bar{\xi }\right) \diag\left( \mu - \bar{
	\xi}\right) \left[I - \left( C-\kappa I\right) Z^{-1}\left( \bar{\xi }\right) \diag\left( \bar{\xi } \right) \right] F \notag \\
	&=& F^{\T} Z^{-1}\left( \bar{\xi }\right) \diag\left( \mu - \bar{
		\xi}\right) \frac{\kappa}{n} \left\{ Z^{\T}\left(\bar{\xi}\right) \right\}^{-1} F \notag \\
	&=& \frac{\kappa}{n} F^{\T} T^{\T} \left( \bar{\xi }\right) \diag\left( \mu - \bar{
		\xi}\right)  T\left(\bar{\xi}\right) F,  \label{eq:dL_dalpha} 
\end{eqnarray}
which is well defined and continuous on $\Xi_{+}$. The equality $$\frac{\kappa}{n} \left\{ Z^{\T}\left(\bar{\xi}\right) \right\}^{-1} = I - \left( C-\kappa I\right) Z^{-1}\left( \bar{\xi }\right) \diag\left( \bar{\xi } \right)$$ for $Z^{\T}\left(\bar{\xi}\right) = \left(C-\kappa I\right) \diag\left( \bar{\xi} \right) +\frac{\kappa}{n} I$ is obtained by applying the identity $(UV + Y)^{-1} = Y^{-1} - Y^{-1} U \left( V Y^{-1} U + I \right)^{-1} V Y^{-1}$ with matrices $U = C - \kappa I$, $V = \diag\left(\bar{\xi}\right)$, and $Y = \frac{\kappa}{n} I$. From the definition of $T(\xi)$ in Section~\ref{subsec:equi_theorem} 
it immediately follows that $\left\{ Z^{\T} \left(\bar{\xi}\right) \right\}^{-1} = T\left(\bar{\xi}\right)$ and $Z^{-1}\left(\bar{\xi}\right) = T^{\T}\left(\bar{\xi}\right)$.

Suppose now that $\bar{\xi} \in \Xi\backslash\Xi_{+}.$ To obtain \eqref{eq:dL_dalpha} from \eqref{eq:L_xi}, we used two standard rules for the derivatives of matrix functions $A(\alpha)$ and $B(\alpha)$, namely $\displaystyle \frac{\partial}{\partial \alpha} \left[ A(\alpha) B(\alpha) \right] = \frac{\partial A(\alpha)}{\partial \alpha} B(\alpha) + A(\alpha) \frac{\partial B(\alpha)}{\partial \alpha}$ and $\displaystyle \frac{\partial A^{-1}(\alpha)}{\partial \alpha} = - A^{-1}(\alpha) \frac{\partial A(\alpha)}{\partial \alpha} A^{-1}(\alpha)$. The same rules hold if instead of the symbol of the derivative $\frac{\partial}{\partial \alpha}$ we use the symbol of the limit $\underset{\alpha \rightarrow 0^+}{\lim}$, abbreviated below by the symbol $\Lambda$. Indeed, 
\begin{eqnarray*}
	\Lambda \left[ A\left( \alpha \right) B\left( \alpha \right) \right]
	&=&\lim_{\alpha \rightarrow 0^{+}}\left[ \frac{A\left( \alpha \right)
		-A\left( 0\right) }\alpha B\left( \alpha \right) +A\left( 0\right) \frac{
		B\left( \alpha \right) -B\left( 0\right) }\alpha \right] \\
	&=&\Lambda \left[ A\left( \alpha \right) \right] B\left( 0\right) +A\left(
	0\right) \Lambda \left[ B\left( \alpha \right) \right],
	\,\,\,\,\,\,\,\,\,\,\,\,\,\,\,\,\,\,\,\,\,\,\,\,\,\,\,\,\,\,\,\,\,\,\,\,\,\,
	\,\,\,\,\,\,\,\,\,\,\,\,\,\,\,\,\,
\end{eqnarray*}
and from $0 = \Lambda\left[A(\alpha) A^{-1}(\alpha) \right] = A(0) \Lambda \left[ A^{-1}(\alpha) \right] + \Lambda \left[A(\alpha) \right] A^{-1}(0)$ we also have
\[
\Lambda \left[ A^{-1}\left( \alpha \right) \right] =-A^{-1}\left( 0\right)
\Lambda \left[ A\left( \alpha \right) \right] A^{-1}\left( 0\right).
\]
Evidently $\Lambda\left[ \diag \left\{ (1-\alpha) \bar{\xi} + \alpha \mu \right\}\right] =\diag\left( \mu - \bar{\xi }\right)$. Hence we finally obtain
that 
\begin{equation}
\lim_{\alpha \rightarrow 0^{+}} \frac{M\left\{ \left( 1-\alpha \right) 
	\bar{\xi}+\alpha \mu \right\} -M\left( \bar{\xi}\right) } {\alpha } =
	\lim_{\alpha \rightarrow 0^{+}} \frac{L\left\{ \left( 1-\alpha \right) 
	\bar{\xi}+\alpha \mu \right\} -L\left( \bar{\xi}\right) } {\alpha } \label{eq:lim_M_alpha0}
\end{equation}
is expressed again by the right-hand side of \eqref{eq:dL_dalpha}, which is evidently a
continuous function of $\bar{\xi }$ on the whole set $\Xi$. \qed

\vspace*{2ex}

\section{More examples and additional results}

%%%%%%%%%%%%%%%%%%%%%%%%%%%%%%%%%%%%%%%%%

	\subsection*{Random sampling results for Examples 1 and 2}

	Tables~\ref{tab:example1-1}, \ref{tab:example1-2}, and \ref{tab:example2} contain further results for Example~1, the modified Example~1, and Example~2 
	when using the random sampling approaches to obtain exact designs. In addition to the approaches introduced in Section~\ref{sec:examples}, we will also consider exact designs found by random sampling from the measure provided by \cite{dette_optimal_2016-1} and denote this approach by R-DET+EP.

\begin{table}[hbtp!]
	\begin{center}
	\caption{Further optimal designs  and D-efficiencies  for Example~1}
	\begin{tabular}{lccccc}
			& $x_1$ & $x_2$ & $x_3$ & $x_4$ & D-eff \\ 
			R-UNIF (median efficiency) &  &  &  &  & 0.6955 \\ 
			R-UNIF (highest efficiency) & 1.12 & 1.30 & 1.72 & 1.96 & 0.8797 \\ 
			R-VN (median efficiency) &  &  &  &  & 0.7746 \\ 
			R-VN (highest efficiency) & 1.23 & 1.69 & 1.79 & 2.00 & 0.9105 \\ 
			R-DPZ+EP (median efficiency) &  &  &  &  & 0.5664 \\ 
			R-DPZ+EP (highest efficiency) & 1.00 & 1.23 & 1.71 & 2.00 & 0.8813 \\  
	\end{tabular}
	\label{tab:example1-1}
	\end{center}
\end{table}

\begin{table}[hbtp!]
	\begin{center}
	\caption{Further optimal designs and D-efficiencies for the modified Example~1}
	\begin{tabular}{lccccc}
			& $x_1$ & $x_2$ & $x_3$ & $x_4$ & D-eff \\ 
			R-UNIF (median efficiency) &  &  &  &  & 0.4887 \\ 
			R-UNIF (highest efficiency) & 1.05 & 1.24 & 1.70 & 1.99 & 0.9207 \\ 
			R-VN (median efficiency) &  &  &  &  & 0.4933 \\ 
			R-VN (highest efficiency) & 1.00 & 1.39 & 1.75 & 2.00 & 0.8405 \\   
	\end{tabular}	\label{tab:example1-2}
	\end{center}
\end{table}

\begin{table}[ht]
	\begin{center}
	\caption{Further optimal designs and D-efficiencies for Example~2}
	\begin{tabular}{lcccccc}
			& $x_1$ & $x_2$ & $x_3$ & $x_4$ & $x_5$ & D-eff \\ 
			R-UNIF (median efficiency) &  &  &  &  &  & 0.3208 \\ 
			R-UNIF (highest efficiency) & 1.04 & 1.11 & 1.28 & 1.80 & 2.00 & 0.8283 \\ 
			R-VN (median efficiency) &  &  &  &  &  & 0.5836 \\ 
			R-VN (highest efficiency) & 1.00 & 1.16 & 1.36 & 1.80 & 2.00 & 0.9299 \\ 
			R-DPZ+EP (median efficiency) &  &  &  &  &  &  0.7409 \\ 
			R-DPZ+EP (highest efficiency) & 1.00 & 1.17 & 1.47 & 1.78 & 2.00 & 0.9281 \\  
	\end{tabular}
	\label{tab:example2}
	\end{center}
\end{table}

\subsection*{Example 4: absolute exponential kernel}\label{sec:example3}

Similar conclusions hold for the next example taken from Section~4.2 of \cite{dette_new_2017}. This four-parameter model is characterized by
\begin{eqnarray*}
	f^\T(x) & = & \left(\sin x, \cos x, \sin 2x, \cos 2x \right), \\
	\Cov\{\varepsilon(x),\varepsilon( x')\} & = & \exp\left( - \left| x - x' \right| \right) \\
	x & \in & [1,2], \\
	\lambda_{\min}(C) & = & 0.005.
\end{eqnarray*}

As \cite{dette_new_2017}, we consider the A-criterion for this example. 
Since the linear programming algorithm requires the criterion to be positive, we select the criterion to be $$\Phi(M) = \left\{\text{tr}\left(M^{-1}\right)\right\}^{-1},$$ for which the derivative is $$\nabla_M \Phi(M) = M^{-2} \cdot \left\{\text{tr}\left(M^{-1}\right)\right\}^{-2},$$
both of which we can plug into the linear Taylor formula 
to obtain the set of linear constraints. 

The selected design points and A-efficiencies for all methods are given in Table~\ref{tab:example3}.

% latex table generated in R 4.0.2 by xtable 1.8-4 package
% Mon Oct 05 16:21:37 2020
\begin{table}[ht]
	\begin{center}
	\caption{Optimal designs and A-efficiencies for Example~4}
	\begin{tabular}{lcccccc}
			& $x_1$ & $x_2$ & $x_3$ & $x_4$ & $x_5$ & A-eff \\ 
			Q-VN & 1.00 & 1.16 & 1.58 & 1.84 & 2.00 & 0.7980 \\ 
			Q-VN+EP & 1.00 & 1.17 & 1.58 & 1.84 & 2.00 & 0.8050 \\ 
			Q-DPZ+EP & 1.00 & 1.25 & 1.50 & 1.75 & 2.00 & 0.7478 \\ 
			R-UNIF (median efficiency) &  &  &  &  &  & 0.0561 \\ 
			R-UNIF (highest efficiency) & 1.01 & 1.13 & 1.57 & 1.88 & 1.99 & 0.6500 \\ 
			R-VN (median efficiency) &  &  &  &  &  & 0.3033 \\ 
			R-VN (highest efficiency) & 1.00 & 1.12 & 1.24 & 1.82 & 2.00 & 0.8555 \\ 
			R-DPZ+EP (median efficiency) &  &  &  &  &  & 0.4764 \\ 
			R-DPZ+EP (highest efficiency) & 1.00 & 1.12 & 1.30 & 1.82 & 2.00 & 0.8414 \\ 
			BKSF & 1.00 & 1.16 & 1.27 & 1.83 & 2.00 & 0.8382 \\ 
			EXS & 1.00 & 1.20 & 1.76 & 1.89 & 2.00 & 0.8602 \\ 
	\end{tabular}
	\label{tab:example3}
	\end{center}
\end{table}

\begin{figure}[hbtp!]
	\centering
	\begin{tabular}{cc}
		\includegraphics[width=0.45\textwidth]{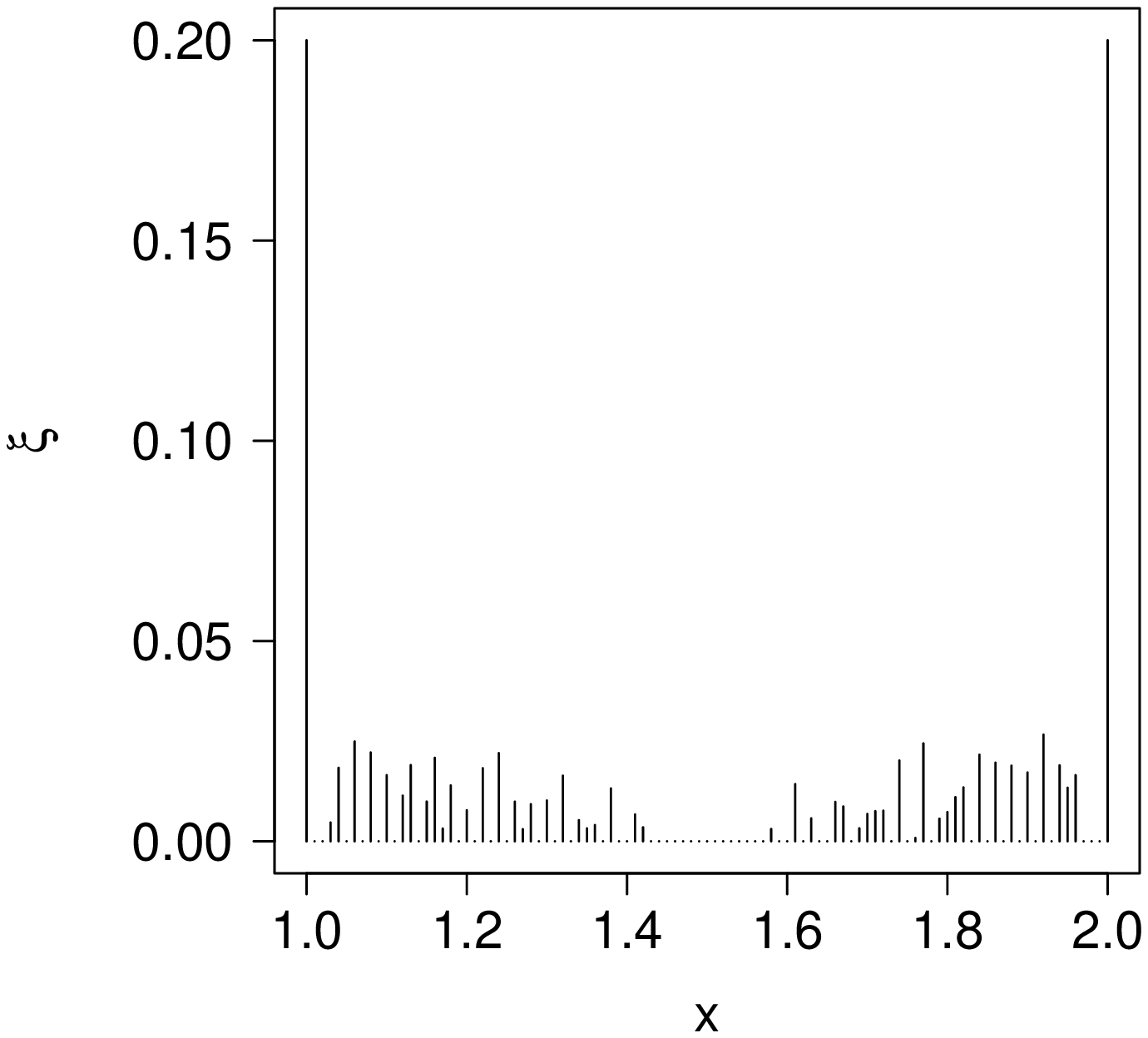} & 
		\includegraphics[width=0.45\textwidth]{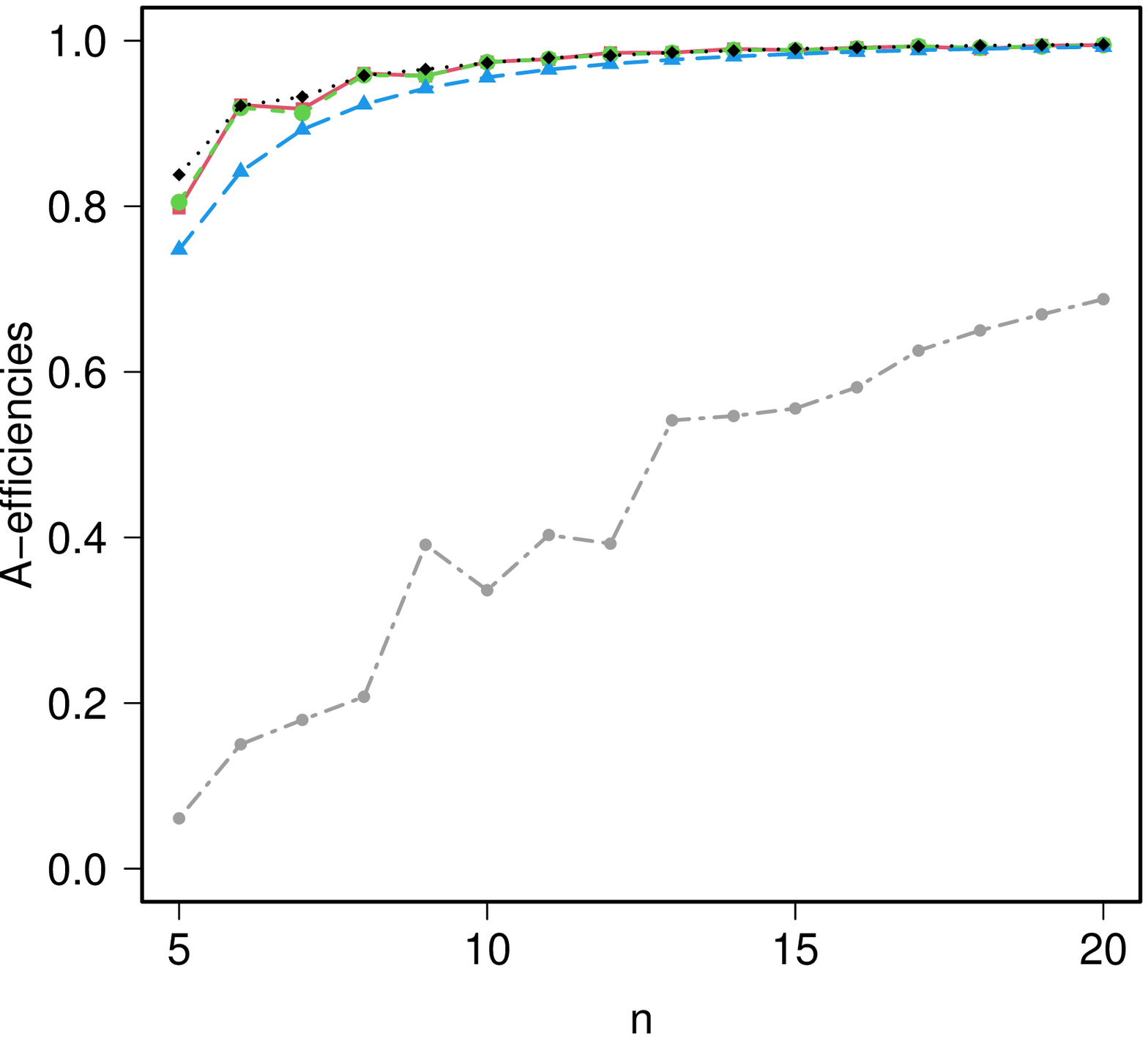}
	\end{tabular}
	\caption{Our measure (left panel) and efficiencies versus sample size (right panel) for the modified Example~4. \label{fig:measures_example3}}
\end{figure}

	\subsection*{Example 5: a bivariate case}\label{sec:example5}

The final example is a multivariate extension of Example~4
to demonstrate once more that our methodology can principally be extended to design dimensions greater than one.

	\begin{eqnarray*}
		x & = & (x_1,x_2)^\T, \\
		f^\T(x) & = & \left(\sin x_1, \cos x_1, \sin 2x_1, \cos 2x_1, \sin x_2, \cos x_2, \sin 2x_2, \cos 2x_2\right), \\
		k(x,x') & = & \exp\left( - \left|x - x' \right| \right), \text{ where }  \left|x - x' \right| = \left|x_1 - x_1'\right| +  \left|x_2 - x_2'\right|\\
		x & \in & [1,2] \times [1,2], \\
		\text{discretized: } x & \in & \{1, 1.1, \ldots, 1.9, 2\} \times \{1, 1.1, \ldots, 1.9, 2\},\\
		\lambda_{\min}(C) & = & 0.002599.
	\end{eqnarray*}

\begin{figure}[hbtp!]
	\centering
	\begin{tabular}{cc}
		\includegraphics[width=0.5\textwidth]{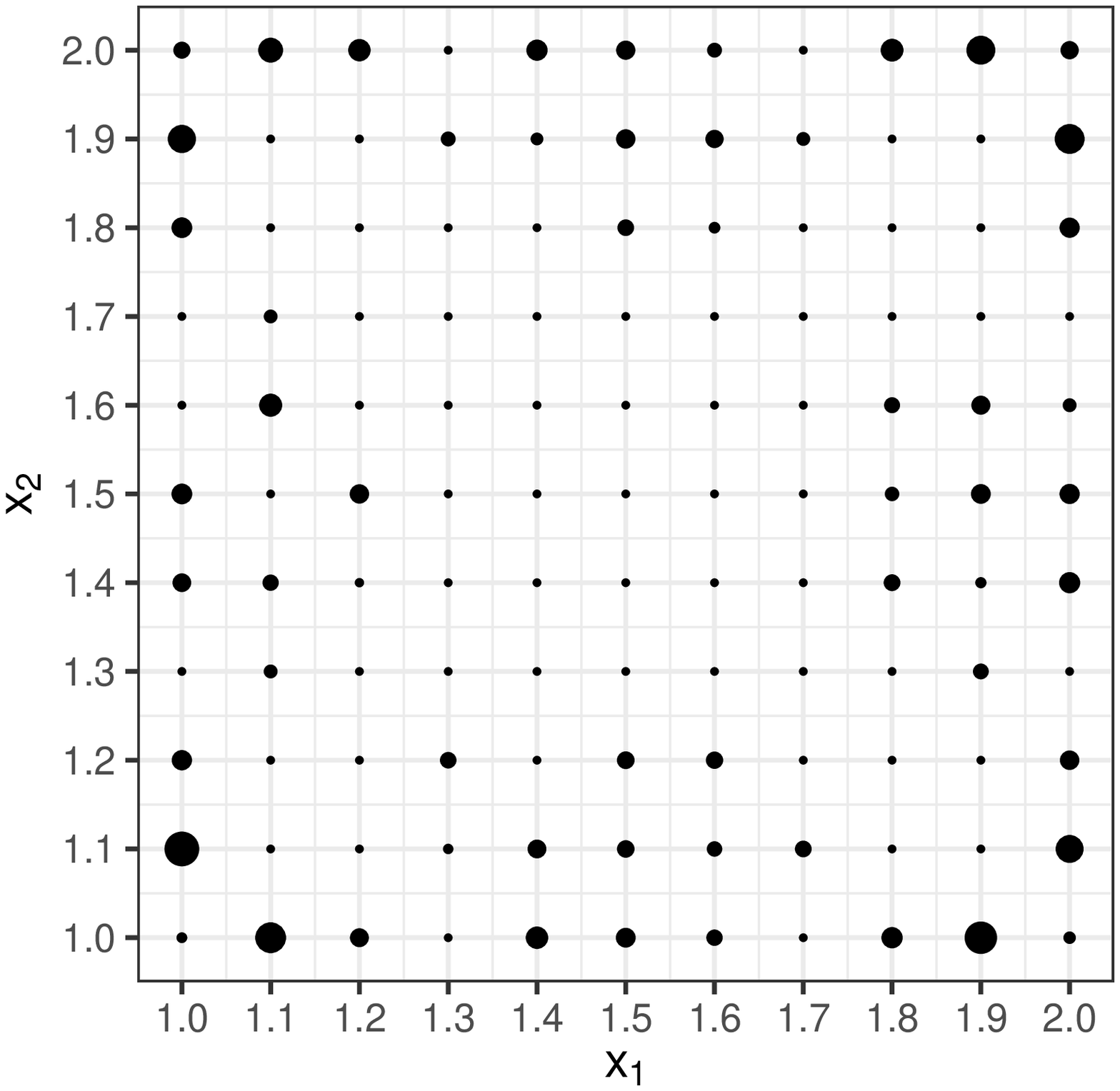} & 
		\includegraphics[width=0.5\textwidth]{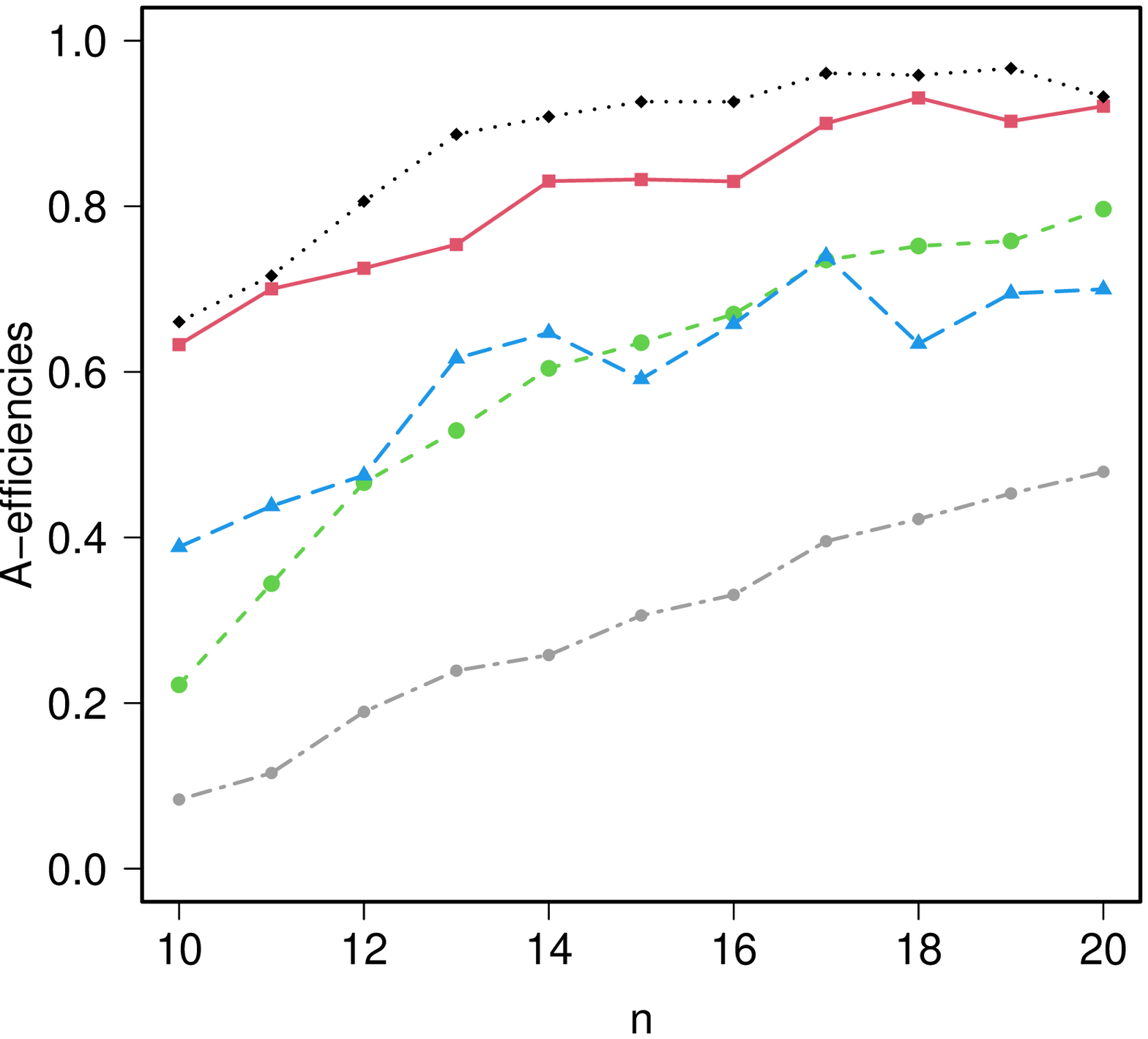}
	\end{tabular}
	\caption{Graphs for Example~5. Left: discrete measure obtained by running linear programming algorithm for virtual noise representation for $n = 10$. Right: A-efficiencies of exact designs with respect to optimal measure of virtual noise representation obtained by various methods for $n = 10$ to $n = 20$. The methods depicted are: R-VN (optimum: solid red line with squares, median: dashed green line with large dots), R-UNIF (optimum: long-dashed blue line with triangles, median: long-short-dashed grey line with small dots), BKSF (dotted black line with diamonds).\label{fig:measures_example5}}
\end{figure}

	To obtain exact designs from our design measure on a two-dimensional grid, we used the random sampling approach. That is, we sampled 100 $n$-point designs according to our measure. In Fig.~\ref{fig:measures_example5}, the best as well as the median A-efficiencies among the sampled designs are plotted. We also sampled 100 designs uniformly on the grid and computed the best and median efficiencies among those designs. Using the best among the sampled designs leads to reasonably efficient designs compared to the algorithm using the approximate sensitivity function proposed by \cite{fedorov_design_1996}.

	\subsection*{Efficiencies for Example~3}

\begin{figure}[hbtp!]
	\centering
	
	\includegraphics[width=0.75\textwidth]{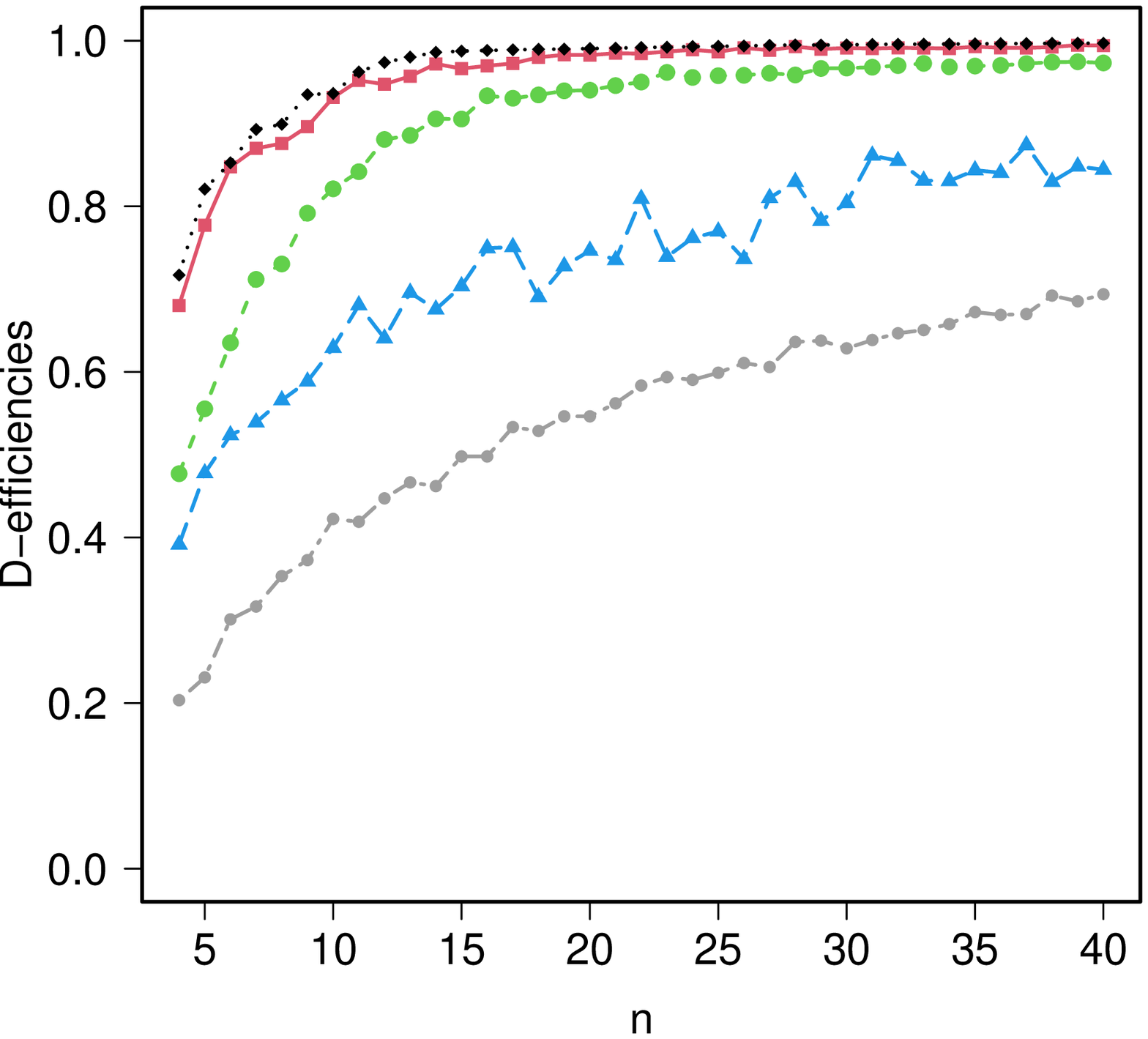} 
	
	\caption{D-efficiencies for Example~3 for $n = 4$ to $n = 40$. The same line types, point symbols and colours are used for the respective methods as in the right panel of Figure~\ref{fig:measures_example5}. \label{fig:efficiencies_example3}}
\end{figure}

	Figure~\ref{fig:efficiencies_example3} displays the D-efficiencies with respect to $\Phi\left\{M \! \left(\bar{\xi}\right)\right\}$ for Example~3 for $n=4$ to $n=40$.

\end{document}